\documentclass[12pt]{article}
\usepackage{mathtools} 
\usepackage{amsmath, amssymb} 
\usepackage{amsthm}
\usepackage{amsrefs}
\usepackage{url} % formats url's nicely using \url{}
\usepackage{enumitem} 

\usepackage{graphicx}
\usepackage{tikz}
\usepackage{tkz-graph}
\usetikzlibrary{shapes}
\usetikzlibrary{arrows}
\usetikzlibrary{decorations.markings}
\usepackage{float}

\usepackage{graphicx}
\usepackage{caption,subcaption}

\newcommand\cx{{\mathbb C}}% complexes

\newcommand\re{{\mathbb R}}%reals

\DeclarePairedDelimiter\abs{\lvert}{\rvert}%
\DeclarePairedDelimiter\norm{\lVert}{\rVert}%

\DeclarePairedDelimiter\bra{\langle}{\rvert}
\DeclarePairedDelimiter\ket{\lvert}{\rangle}

\makeatletter
\let\oldabs\abs
\def\abs{\@ifstar{\oldabs}{\oldabs*}}
\let\oldnorm\norm
\def\norm{\@ifstar{\oldnorm}{\oldnorm*}}
\makeatother

\newcommand\comp[1]{{\mkern2mu\overline{\mkern-2mu#1}}}

\newtheoremstyle{plainsl}%
	{\topsep}
	{\topsep}
	{\slshape} % only non-default setting
	{}
	{\normalfont\bfseries}
	{.}
	{ }
	{}

\swapnumbers

{\theoremstyle{plainsl}
\newtheorem{theorem}{Theorem}[section]
\newtheorem{lemma}[theorem]{Lemma}
\newtheorem{corollary}[theorem]{Corollary}}
{\theoremstyle{remark}
}

\renewcommand\proof{\noindent\textsl{Proof. }}
\newcommand\sqr[2]{{\vbox{\hrule height.#2pt
    \hbox{\vrule width.#2pt height#1pt \kern#1pt
        \vrule width.#2pt}\hrule height.#2pt}}}
\renewcommand\qed{%
	\ifmmode\eqno\sqr53
	\else\nolinebreak\ \hfill\sqr53\medbreak\fi}

\DeclareMathOperator{\rk}{rk}
\DeclareMathOperator{\tr}{tr}
\DeclareMathOperator{\col}{Col}

\newcommand\aut[1]{{\rm Aut}(#1)}
\usepackage{blkarray}

\usepackage{hyperref}

\title{Pair State Transfer}
\author{Qiuting Chen, Chris Godsil\footnote{Research supported by Natural Sciences and Engineering 
	Council of Canada, Grant No. RGPIN-9439}\\
 	Department of Combinatorics \& Optimization,\\ 
	University of Waterloo, Waterloo, Ontario, Canada}

\begin{document}
\maketitle
\abstract{Let $L$ denote the Laplacian matrix of a graph $G$. We study continuous quantum walks on $G$ defined by the transition matrix $U(t)=\exp\left(itL\right)$. The initial state is of the pair state form, $e_a-e_b$ with $a,b$ being any two vertices of $G$. We provide two ways to construct infinite families of graphs that have perfect pair state transfer. We study a ``transitivity" phenomenon which cannot occur in vertex state transfer. We characterize perfect pair state transfer on paths and cycles. We also study the case when quantum walks are generated by the unsigned Laplacians of underlying graphs and the initial states are of the plus state form, $e_a+e_b$. When the underlying graphs are bipartite, plus state transfer is equivalent to pair state transfer. }

\section{Introduction}

The concept of a quantum walk was introduced by Farhi and Gutmann \cite{Farhi1998}  as a quantum mechanical analogue of a classical random walk on decision trees. Exploiting the interference effects of quantum mechanics, quantum walks outperform classical random walks for some computational tasks~\cite{Childs2002}.

In the field of quantum information processing, Christandl et al.~\cite{Christandl2004} brought our attention to the topic of perfect state transfer.   Using the tool of quantum scattering theory, Childs~\cite{Childs2008} proved that continuous time quantum walks can be regarded as a universal computational primitive and any desired quantum computation can be encoded in some underlying graph of the quantum walk. Quantum walks have become powerful tools to improve existing quantum algorithms and develop new quantum algorithms. In this paper, we use graphs to represent networks of interacting qubits and study quantum state transfer during quantum communication over the network.

Let $G$ be a graph. The evolution of a continuous quantum walk on $G$ is given
by the matrices
\[
	U(t) = \exp(itH), \qquad (t\in \re). 
\]
Here $H$ is a matrix, called the Hamiltonian of the walk, and is usually either the
adjacency matrix, the Laplacian, or the signless Laplacian of $G$. In any case, $H$ is Hermitian
and its rows and columns are indexed by the vertices of $G$. If $n=|V(G)|$, then the walk
represents a quantum system with an $n$-dimensional state space. We identify the states of the system
by density matrices, i.e., positive semidefinite matrices with trace 1.
The physically meaningful questions are the form: ``Given that the quantum system is
initially in state represented by a density matrix $D_0$, what is the probability that, at time $t$,
its state is $D_1$?''.

Let $e_1,e_2,\cdots,e_n$ denote the standard basis for $\cx^n$. Then $D_r=e_re_r^T$ is a density matrix,
and the mathematical questions reduce to question about the absolute value of $U(t)_{r,s}$
(for given vertices $r$ and $s$). Thus if $r\ne s$ and $|U(t)_{r,s}|=1$ at some time $t$,
we say that it has perfect state transfer from vertex $r$ to vertex $s$ and this is equivalent to having perfect state transfer from vertex $s$ to vertex $r$ since $U(t)$ is symmetric. Perfect state transfer
is a potentially useful tool in quantum computation, and so there is a considerable literature
on the topic \cite{Bose2003,Childs2008,Childs2013}. In \cite{Godsil2017}, it was observed that 
most of the results on the topic only require the fact that the density matrices $e_re_r^T$ are real. 
This suggests strongly that there may be other density matrices of interest. In this paper we focus 
on density matrices of the form
\[
	(e_a-e_b)(e_a-e_b)^T,
\]
for vertices $a$ and $b$ in some graph; we call this a \textsl{pair state}. Such a density matrix 
is the Laplacian matrix for a graph formed from a single edge with $a,b$ being its ends and hence it 
seems natural to consider continuous walks using graph Laplacians as Hamiltonians.
The goal of this paper is to investigate the properties of these walks.

We prove that perfect pair state transfer is preserved under complementations and under taking Cartesian 
product. This helps us to construct more examples of graphs with perfect pair state transfer.

Transitivity is one phenomenon that only can occur when the initial state is a pair state and this
will be discussed in Section~\ref{transitivity}.

Although pair state transfer has monogamy 
and symmetry properties just as vertex state transfer does, because our initial state involves two 
vertices, rather than just one, we have more flexibility. One consequence is that we have 
more examples of perfect state transfer using pair states as the initial state. The following 
table computationally indicates that perfect pair state transfer occurs more often than perfect 
vertex state transfer in the same set of graphs.
\begin{table}[H]
\begin{center}
\begin{tabular}{|c|c|c|c|c|c|}
\hline
$G_n$ & Total  & vertex PST 
 &Prop. &  pair PST 
 &Prop.\\
\hline
$G_5$ & $21$  & $1$ & $4.8\%$&$6$ & $28.6\%$ \\
\hline
$G_6$ & $112$ &$1$&$0.9\%$ & $27$ & $24.1\%$\\
\hline
$G_7$ & $853$  &$1$& $0.1\%$ &$104$ & $12.2\%$\\
\hline
$G_8$ & $11117$ & $ 5$ & $0.004\%$ & $779$ & $7.0\%$\\
\hline
\end{tabular}
\end{center}
\caption{the number of graphs with vertex PST versus the number of graphs with pair PST}
\label{periodic vertices}
\end{table}

Let $G_n$ denote the set of connected graphs on $n$ vertices and the second column shows the cardinality of the corresponding $G_n$. The third column shows the number of graphs that have perfect vertex state transfer in $G_n$ when the adjacency matrices are the Hamiltonians and the fourth column shows the corresponding proportion. The fifth column shows the number of graphs that have perfect pair state transfer in $G_n$ when the Laplacians are the Hamiltonians and the last column shows the corresponding proportion. In Section~\ref{complement}, we prove that perfect pair state transfer is preserved under taking the complement of the underlying graphs. So to avoid overcounting, the tables in this paper only show the number of graphs with perfect pair state transfer between pairs such that at least one of them is an edge.

Perfect state transfer is a significant phenomenon in quantum communication, but quite rare in quantum walks. We always want to find more graphs with perfect state transfer. On a fixed number of vertices, there are more graphs with perfect pair state transfer than graphs with perfect vertex state transfer. This is a huge advantage of Laplacian pair state transfer and this is also why Laplacian pair state transfer is interesting.

We prove that perfect edge state transfer occurs on $P_n$ if and only if $n=3$ or $4$.  We also prove 
that $C_4$ is the only cycle that has perfect edge state transfer. 

We also study the case when the unsigned Laplacian of the underlying graph is used as Hamiltonian and 
the initial state of the form $e_a+e_b$. We prove that perfect pair state transfer is equivalent to 
perfect plus state transfer when the underlying graph is bipartite. This allows us to prove analogous results 
for perfect plus state transfer. That is, $P_3$ and $P_4$ are the only paths and $C_4$ is the only 
cycle with perfect plus state transfer.

\section{Preliminaries}
\subsection{Pair State Transfer}
Let $G$ be a graph with $n$ vertices. 
Let $A$ denote the adjacency matrix of $G$ and let $\Delta$ denote the degree matrix of $G$ . Then the Laplacian of $G$ is the matrix such that \[
L=\Delta-A.\] When $L$ is the Hamiltonian associated to the quantum walk on $G$, the transition matrix is \[
U(t) = \exp\left(itL\right).\] By Schr\"{o}dinger’s equation, the probability of the state at $\ket{a_2}$ starting at $\ket{a_1}$ after time $t$ is 
\[\abs{ \bra{ a_2}U(t)\ket{a_1}}^2.\]

The quantum pair state associated with a pair of vertices  $(a,b)$ of $G$ is represented by \[e_a-e_b,\] where $e_a,e_b\in \mathbb{R}^n$ are the characteristic vectors of $a,b$ respectively. We want quantum states to be represented by unit vectors, so when we perform computations about pair state transfer, we use the normalized pair state \[
\frac{1}{\sqrt{2}}(e_a-e_b),\] but except for that, in this paper we always use $e_a-e_b$ to denote our pair states for convenience. Unless explicitly stated otherwise, the initial state is a pair state in the continuous quantum walk on $G$ generated by the Laplacian of $G$.

There is perfect pair state transfer between $(a,b)$ and $(c,d)$ if there exists a complex scalar $\gamma$ with $\abs{\gamma}=1$ satisfying that 
\[
	U(t)(e_a-e_b)=\gamma(e_c-e_d)
\] for some non-negative time $t$ and probabilistically 
\[
	\abs{\frac{1}{2}(e_c-e_d)^TU(t)(e_a-e_b)}^2=1.
\] 
We say the pair $(a,b)$ is periodic with period $\tau$ if it has perfect state transfer to itself at time $\tau$. 

Another way to represent a quantum state is density matrices. A density matrix is a positive semidefinite matrix of trace $1$. A density matrix $D$ represents a pure state if $\rk(D)=1$ or equivalently $\tr(D^2)=1$ . Let $e_i$ denote the standard basis vector in $\mathbb{C}^{\abs{V(G)}}$ indexed by the vertex $i$ in graph $G$ and then \[D=\frac{1}{2}(e_a-e_b)(e_a-e_b)^T\] is  a pure state associated with a pair of vertices $(a,b)$ in $G$, which we call the density matrix of pair $(a,b)$.

Given a density matrix $D$ as the initial state of a continuous quantum walk, then the state that $D$ is transferred to at time $t$ is given by 
\[ D(t) = U(t)DU(-t),\] 
where $U(t)=\exp(itL)$ is the usual transition matrix associated with graph $G$ whose Laplacian matrix is $L$. There is perfect state transfer between density matrices $P$ and $Q$, which means that there is a time $t$ such that 
\[
Q = U(t)PU(-t).\] 
We say a state $P$ is periodic if there is a time $t$ such that 
\[
P= U(t)PU(-t).\]

\begin{figure}[H]
\begin{center}
\begin{tikzpicture}[rotate=90]
\definecolor{cv0}{rgb}{0.0,0.0,0.0}
\definecolor{cfv0}{rgb}{1.0,1.0,1.0}
\definecolor{clv0}{rgb}{0.0,0.0,0.0}
\definecolor{cv1}{rgb}{0.0,0.0,0.0}
\definecolor{cfv1}{rgb}{1.0,1.0,1.0}
\definecolor{clv1}{rgb}{0.0,0.0,0.0}
\definecolor{cv2}{rgb}{0.0,0.0,0.0}
\definecolor{cfv2}{rgb}{1.0,1.0,1.0}
\definecolor{clv2}{rgb}{0.0,0.0,0.0}
\definecolor{cv3}{rgb}{0.0,0.0,0.0}
\definecolor{cfv3}{rgb}{1.0,1.0,1.0}
\definecolor{clv3}{rgb}{0.0,0.0,0.0}
\definecolor{cv4}{rgb}{0.0,0.0,0.0}
\definecolor{cfv4}{rgb}{1.0,1.0,1.0}
\definecolor{clv4}{rgb}{0.0,0.0,0.0}
\definecolor{cv5}{rgb}{0.0,0.0,0.0}
\definecolor{cfv5}{rgb}{1.0,1.0,1.0}
\definecolor{clv5}{rgb}{0.0,0.0,0.0}
\definecolor{cv0v3}{rgb}{0.0,0.0,0.0}
\definecolor{cv0v4}{rgb}{0.0,0.0,0.0}
\definecolor{cv1v4}{rgb}{0.0,0.0,0.0}
\definecolor{cv1v5}{rgb}{0.0,0.0,0.0}
\definecolor{cv2v4}{rgb}{0.0,0.0,0.0}
\definecolor{cv2v5}{rgb}{0.0,0.0,0.0}
\definecolor{cv3v5}{rgb}{0.0,0.0,0.0}
\Vertex[style={minimum size=1.0cm,draw=cv0,fill=cfv0,text=clv0,shape=circle},LabelOut=false,L=\hbox{$0$},x=0cm,y=6.0cm]{v0}
\Vertex[style={minimum size=1.0cm,draw=cv1,fill=cfv1,text=clv1,shape=circle},LabelOut=false,L=\hbox{$1$},x=0.0cm,y=0cm]{v1}
\Vertex[style={minimum size=1.0cm,draw=cv2,fill=cfv2,text=clv2,shape=circle},LabelOut=false,L=\hbox{$2$},x=3cm,y=0.0cm]{v2}
\Vertex[style={minimum size=1.0cm,draw=cv3,fill=cfv3,text=clv3,shape=circle},LabelOut=false,L=\hbox{$3$},x=3.0cm,y=6cm]{v3}
\Vertex[style={minimum size=1.0cm,draw=cv4,fill=cfv4,text=clv4,shape=circle},LabelOut=false,L=\hbox{$4$},x=0cm,y=3cm]{v4}
\Vertex[style={minimum size=1.0cm,draw=cv5,fill=cfv5,text=clv5,shape=circle},LabelOut=false,L=\hbox{$5$},x=3cm,y=3cm]{v5}
\Edge[lw=0.01cm,style={color=cv0v3,},](v0)(v3)
\Edge[lw=0.01cm,style={color=cv0v4,},](v0)(v4)
\Edge[lw=0.01cm,style={color=cv1v4,},](v1)(v4)
\Edge[lw=0.01cm,style={color=cv1v5,},](v1)(v5)
\Edge[lw=0.01cm,style={color=cv2v4,},](v2)(v4)
\Edge[lw=0.01cm,style={color=cv2v5,},](v2)(v5)
\Edge[lw=0.01cm,style={color=cv3v5,},](v3)(v5)
\end{tikzpicture}
\caption{Smallest graph with PST from an edge pair to a non-edge pair}
\label{pair pst}
\end{center}
\end{figure}
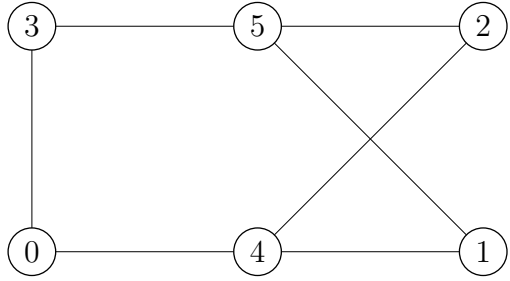

Here, we want to emphasize that the pair of vertices $(a,b)$ associated with the pair state $e_a-e_b$, need not be adjacent. We say that $e_a-e_b$ is an edge state only when $(a,b)$ is an edge. There may perfect pair state transfer between $(a,b)$ and $(c,d)$ where $(a,b)$ is an edge while $(c,d)$ is not.
Thus there is perfect pair state transfer between $(0,3)$ and $(4,5)$ in the graph $G$ shown in Figure~\ref{pair pst}.

The laplacian matrix of $G$ is 
\[ L=
\begin{pmatrix}
2 & 0 & 0 & -1 & -1 & 0 \\
0 & 2 & 0 & 0 & -1 & -1 \\
0 & 0 & 2 & 0 & -1 & -1 \\
-1 & 0 & 0 & 2 & 0 & -1 \\
-1 & -1 & -1 & 0 & 3 & 0 \\
0 & -1 & -1 & -1 & 0 & 3
\end{pmatrix}\] 
and it has Laplacian eigenvalues $\{3\pm \sqrt{3},0,2,4\}$. The pair state $e_0-e_3$ and $e_4-e_5$ both have eigenvalue support $\{2,4\}$. Then   
\[
	(e_0-e_3)^T\exp(itL)(e_4-e_5)
	 	=e^{2it}(e_0-e_3)^TE_2(e_4-e_5)+e^{4it}(e_0-e_3)^TE_4(e_4-e_5)
\]
and since
\[
E_2=\begin{pmatrix}
\frac{1}{4} & 0 & 0 & -\frac{1}{4} & \frac{1}{4} & -\frac{1}{4} \\
0 & \frac{1}{2} & -\frac{1}{2} & 0 & 0 & 0 \\
0 & -\frac{1}{2} & \frac{1}{2} & 0 & 0 & 0 \\
-\frac{1}{4} & 0 & 0 & \frac{1}{4} & -\frac{1}{4} & \frac{1}{4} \\
\frac{1}{4} & 0 & 0 & -\frac{1}{4} & \frac{1}{4} & -\frac{1}{4} \\
-\frac{1}{4} & 0 & 0 & \frac{1}{4} & -\frac{1}{4} & \frac{1}{4}
\end{pmatrix},\quad 
E_4=\begin{pmatrix}
\frac{1}{4} & 0 & 0 & -\frac{1}{4} & -\frac{1}{4} & \frac{1}{4} \\
0 & 0 & 0 & 0 & 0 & 0 \\
0 & 0 & 0 & 0 & 0 & 0 \\
-\frac{1}{4} & 0 & 0 & \frac{1}{4} & \frac{1}{4} & -\frac{1}{4} \\
-\frac{1}{4} & 0 & 0 & \frac{1}{4} & \frac{1}{4} & -\frac{1}{4} \\
\frac{1}{4} & 0 & 0 & -\frac{1}{4} & -\frac{1}{4} & \frac{1}{4}
\end{pmatrix},
\] 
we get that 
\begin{align*}
(e_0-e_3)^TU(t)(e_4-e_5)
&=e^{2it}\left(\frac{1}{4}+\frac{1}{4}+\frac{1}{4}+\frac{1}{4}\right)+e^{4it}\left(-\frac{1}{4}-\frac{1}{4}-\frac{1}{4}-\frac{1}{4}\right)\\
&=e^{2it}-e^{4it}.
\end{align*} 
So at time $t=\frac{\pi}{2}$, we have  
\[
\abs{\frac{1}{2}(e_0-e_3)^TU(\frac{\pi}{2})(e_4-e_5)}^2=\abs{\frac{1}{2}\left(
e^{\pi i}-e^{2\pi i}\right)}^2=1.
\] 
This shows that when $t=\frac{\pi}{2}$, there is perfect state transfer between $e_0-e_3$ and $e_4-e_5$.

\subsection{Transition Matrices}
\label{trans mtx}
Let $G$ be a graph with the Laplacian matrix $L$. We assume that the eigenvalues of $L$ are $\theta_1,\theta_2,\cdots,\theta_n$. We know that $L$ is a real symmetric matrix.
Then the spectral decomposition of $L$ is  \begin{equation}
\label{spd}
L=\sum_{i=1}^n \theta_i E_i,
\end{equation}where the matrices $E_1,E_2,\cdots,E_n$ satisfy:\begin{enumerate}[label=(\roman*)]
\item $\sum_{i=1}^n E_i = I$,
\item $E_rE_s=\begin{cases}
0, \quad\text{if } r\neq s;\\
E_r, \quad\text{if } r=s.
\end{cases}$
\end{enumerate}
A matrix $E$ is an idempotent if $E^2=E$. The matrices $E_1,E_2,\cdots,E_n$ in Equation~\ref{spd} are called spectral idempotents and $E_r$ represents the orthogonal projection onto the $\theta_r$-eigenspace of $L$. 
\begin{theorem}
\label{spd-function}
Let $M$ be a real symmetric matrix and let $\sum_{i=1}^n\theta_i E_i$ denote the spectral decomposition of $M$. If $f(x)$ is a function defined on the eigenvalues of $M$, then 
\[
f(M)=\sum_{i=1}^n f(\theta_i)E_i.\qed
\]
\end{theorem} 

This standard algebraic graph theory result helps us to obtain the spectral decomposition of the transition matrix $U(t)$, which brings the spectrum of the underlying graph into the picture of continuous quantum walk.

Let $\sum_{i=1}^n \theta_iE_i$ denote the spectral decomposition of the Laplacian matrix of a graph $G$. Then the transition matrix of pair state transfer on $G$ is 
\begin{equation}
\label{spd_U(t)}
U(t) = \sum_{i=1}^n e^{it\theta_i} E_i.
\end{equation}

For example, the spectral decomposition of the Laplacian $L$ of $P_3$ is 
\begin{equation*}
L=0
\begingroup
\renewcommand*{\arraystretch}{1.5}
\begin{pmatrix}
\frac{1}{3}&\frac{1}{3}&\frac{1}{3}\\
\frac{1}{3}&\frac{1}{3}&\frac{1}{3}\\
\frac{1}{3}&\frac{1}{3}&\frac{1}{3}
\end{pmatrix}
\endgroup
+
1\begingroup
\renewcommand*{\arraystretch}{1.5}
\begin{pmatrix}
\frac{1}{2}&0&-\frac{1}{2}\\
0&0&0\\
-\frac{1}{2}&0&\frac{1}{2}
\end{pmatrix}
\endgroup+3
\begingroup
\renewcommand*{\arraystretch}{1.5}
\begin{pmatrix}
\frac{1}{6}&-\frac{1}{3}&\frac{1}{6}\\
-\frac{1}{3}&\frac{2}{3}&-\frac{1}{3}\\
\frac{1}{6}&-\frac{1}{3}&\frac{1}{6}
\end{pmatrix}.
\endgroup
\end{equation*} 
Then the transition matrix associated is\begin{align*}
U(t)&=e^{0it}
\begingroup
\renewcommand*{\arraystretch}{1.5}
\begin{pmatrix}
\frac{1}{3}&\frac{1}{3}&\frac{1}{3}\\
\frac{1}{3}&\frac{1}{3}&\frac{1}{3}\\
\frac{1}{3}&\frac{1}{3}&\frac{1}{3}
\end{pmatrix}
\endgroup+
e^{it}
\begingroup
\renewcommand*{\arraystretch}{1.5}
\begin{pmatrix}
\frac{1}{2}&0&-\frac{1}{2}\\
0&0&0\\
-\frac{1}{2}&0&\frac{1}{2}
\end{pmatrix}
\endgroup+e^{3it}
\begingroup
\renewcommand*{\arraystretch}{1.5}
\begin{pmatrix}
\frac{1}{6}&-\frac{1}{3}&\frac{1}{6}\\
-\frac{1}{3}&\frac{2}{3}&-\frac{1}{3}\\
\frac{1}{6}&-\frac{1}{3}&\frac{1}{6}
\end{pmatrix}
\endgroup \\[1em]
&=
\begingroup
\renewcommand*{\arraystretch}{1.5}
\begin{pmatrix}
\frac{1}{3}+\frac{1}{2}e^{it}+\frac{1}{6}e^{3it}&\frac{1}{3}-\frac{1}{3}e^{3it}&\frac{1}{3}-\frac{1}{2}e^{it}+\frac{1}{6}e^{3it}\\
\frac{1}{3}-\frac{1}{3}e^{3it}&\frac{1}{3}+\frac{2}{3}e^{3it}&\frac{1}{3}-\frac{1}{3}e^{3it}\\
\frac{1}{3}-\frac{1}{2}e^{it}+\frac{1}{6}e^{3it}&\frac{1}{3}-\frac{1}{3}e^{3it}&\frac{1}{3}+\frac{1}{2}e^{it}+\frac{1}{6}e^{3it}
\end{pmatrix}
\endgroup
\end{align*}
\begin{figure}[H]
\begin{center}
\begin{tikzpicture}[scale=0.9]
\definecolor{cv0}{rgb}{0.0,0.0,0.0}
\definecolor{cfv0}{rgb}{1.0,1.0,1.0}
\definecolor{clv0}{rgb}{0.0,0.0,0.0}
\definecolor{cv1}{rgb}{0.0,0.0,0.0}
\definecolor{cfv1}{rgb}{1.0,1.0,1.0}
\definecolor{clv1}{rgb}{0.0,0.0,0.0}
\definecolor{cv2}{rgb}{0.0,0.0,0.0}
\definecolor{cfv2}{rgb}{1.0,1.0,1.0}
\definecolor{clv2}{rgb}{0.0,0.0,0.0}
\definecolor{cv0v1}{rgb}{0.0,0.0,0.0}
\definecolor{cv1v2}{rgb}{0.0,0.0,0.0}
\Vertex[style={minimum size=1.0cm,draw=cv0,fill=cfv0,text=clv0,shape=circle},LabelOut=false,L=\hbox{$0$},x=0.0cm,y=0cm]{v0}
\Vertex[style={minimum size=1.0cm,draw=cv1,fill=cfv1,text=clv1,shape=circle},LabelOut=false,L=\hbox{$1$},x=6cm,y=0cm]{v1}
\Vertex[style={minimum size=1.0cm,draw=cv2,fill=cfv2,text=clv2,shape=circle},LabelOut=false,L=\hbox{$2$},x=12.0cm,y=0cm]{v2}
\Edge[lw=0.01cm,style={color=cv0v1,},](v0)(v1)
\Edge[lw=0.01cm,style={color=cv1v2,},](v1)(v2)
\end{tikzpicture}
\caption{$P_3$}
\end{center}
\end{figure}
When $t=\frac{\pi}{2}$, the transition matrix of $P_3$ is\[U\left(\frac{\pi}{2}\right)=\begingroup
\renewcommand*{\arraystretch}{1.5}
\begin{pmatrix}
\frac{1}{3}+\frac{1}{3}i &\frac{1}{3}+\frac{1}{3}i &\frac{1}{3}-\frac{2}{3}i\\
\frac{1}{3}+\frac{1}{3}i & \frac{1}{3}-\frac{2}{3}i &\frac{1}{3}+\frac{1}{3}i \\
\frac{1}{3}-\frac{2}{3}i &\frac{1}{3}+\frac{1}{3}i &\frac{1}{3}+\frac{1}{3}i 
\end{pmatrix}
\endgroup\]
Then we can see that \begin{align*}
&\abs{\frac{1}{2}(e_1-e_2)U\left(\frac{\pi}{2}\right)(e_0-e_1)}^2\\
&=\abs{\frac{1}{2}\left(\frac{1}{3}+\frac{1}{3}i-\left(\frac{1}{3}-\frac{2}{3}i\right)-\left(\frac{1}{3}-\frac{2}{3}i\right)+\frac{1}{3}+\frac{1}{3}i\right)}^2\\
&=\abs{\frac{1}{2}\left(2i\right)}^2\\
&=1,
\end{align*} which implies that there is perfect pair state transfer from $e_0-e_1$ to $e_1-e_2$ at time $\frac{\pi}{2}$ in $P_3$. Later we will prove actually $P_3$ and $P_4$ are the only paths that have perfect edge state transfer.
\subsection{Eigenvalue Supports}
From Equation~\ref{spd_U(t)}, we also can see that the Laplacian eigenvalues of $G$ play a large role in the pair state transfer. Let $E_r$ be a spectral idempotent such that \[
E_r(e_a-e_b)=0.\] Then we can see that when we talk about the state transfer started in the state $e_a-e_b$, the eigenvalue $\theta_r$ and its idempotent $E_r$ contribute nothing to the evolution. 

The eigenvalue support of the state $e_a-e_b$ is the set of Laplacian eigenvalues $\theta_r$ such that the corresponding idempotent $E_r$ satisfies \[
E_r(e_a-e_b)\neq 0.\]
Thus, when we talk about quantum state transfer initialized in the state $e_a-e_b$, we only care about the eigenvalues in the eigenvalue support of $e_a-e_b$. Recall the example of $P_3$ in previous section. From the spectral decomposition of the Laplacian of $P_3$, we can see that the eigenvalue supports of pair states $e_0-e_1$ and $e_1-e_2$ are the same, that is, $\{1,3\}$.

We say two states $e_a-e_b$ and $e_c-e_d$ are strongly cospectral in $G$ if and only if for each spectral idempotent $E_r$ of the Laplacian of $G$, we have \[
E_r(e_a-e_b)=\pm E_r(e_c-e_d).\] Thus, we can see that if two states are strongly cospectral in graph $G$, then their eigenvalue supports are the same. From the definition of strong cospectrality, the theorem and the corollary below follows immediately. 

\begin{theorem}
\label{pst-strongly cospectral}
If there is perfect state transfer between $e_a-e_b$ and $e_c-e_d$ in graph $G$, then $e_a-e_b$ and $e_c-e_d$ are strongly cospectral.\qed
\end{theorem}

\begin{corollary}
\label{pststronlycospec}
If there is perfect state transfer between $e_a-e_b$ and $e_c-e_d$ in $G$, then $e_a-e_b$ and $e_c-e_d$ have the same eigenvalue support.\qed
\end{corollary}
Now we let $\Lambda_{ab,cd}^+$ denote the set of eigenvalues such that \[
E_r(e_a-e_b)=E_r(e_c-e_d)\] and let $\Lambda_{ab,cd}^-$ denote the set of eigenvalues such that \[
E_r(e_a-e_b)=-E_r(e_c-e_d).\] It is easy to see that \[
\Lambda_{ab}=\Lambda_{cd}=\Lambda_{ab,cd}^+\cup \Lambda_{ab,cd}^-, \quad \Lambda_{ab,cd}^+\cap \Lambda_{ab,cd}^-=\emptyset.\]

Using strong cospectrality, we can derive a characterization of perfect pair state transfer. Since the proof is similar to the proof of perfect vertex state transfer, we omit the proof here. One can refer to \cite[Theorem~2.4.2]{Coutinho2014a} for details.
\begin{lemma}
\label{pst first iff}
Let $X$ be a graph and $a,b,c,d\in V(X)$. Perfect pair state transfer between $e_a-e_b$ and $e_c-e_d$ occurs at time $\tau$ if and only if all of the following conditions hold.
\begin{enumerate}[label=(\alph*)]
\item Pair states $e_a-e_b$ and $e_c-e_d$ are strongly cospectral. Let $\theta_0\in \Lambda_{ab,cd}^+$.
\item For all $\theta_r\in \Lambda_{ab,cd}^+$, there is a $k$ such that $\tau(\theta_0-\theta_r)=2k\pi$.
\item For all $\theta_r\in \Lambda_{ab,cd}^-$, there is a $k$ such that $\tau(\theta_0-\theta_r)=(2k+1)\pi$.\qed
\end{enumerate}
\end{lemma}
\section{Perfect State Transfer \& Periodicity}
In this section, we introduce symmetry and monogamy properties of perfect pair state transfer and give a characterization of periodicity in terms of eigenvalues. We also give a characterization of a fixed state in pair state transfer, which is a special case of periodicity.
\subsection{Basic Properties}
\label{basic properties}
The original results in this section can be found in \cite{Godsil2010} stated and proved in terms of vertex states.
Since proofs of the results using pair states are very similar to the proofs using vertex states, here we just state the results without proofs. One can see \cite{CQTmmath} for detailed proofs using pair states.

Just like perfect vertex state transfer, symmetry and monogamy are two basic properties of perfect pair state transfer.
\begin{theorem}
There is perfect state transfer from $e_a-e_b$ to $e_c-e_d$ in graph $G$ at time $\tau$ if and only if there is perfect state transfer from  $e_c-e_d$  to $e_a-e_b$ at time $\tau$. \qed
\end{theorem}
\begin{theorem}
\label{pst-ped}
Suppose that $e_a-e_b$ has perfect state transfer at time $\tau$ in graph $G$. Then $e_a-e_b$ is periodic at time $2\tau$.\qed
\end{theorem}

Thus, being periodic is a necessary condition for a state to have perfect state transfer.

\begin{theorem}
If there is perfect state transfer between  $(a,b)$ and $(c,d)$ in graph $G$, then both $(a,b)$ and $(c,d)$ are periodic with the same minimum period. If the minimum period is $\sigma$, then perfect state transfer between the two edges occurs at time $\frac{1}{2}\sigma$. \qed
\end{theorem}
Since perfect state transfer occurs exactly at half of the period, the monogamy property of perfect state transfer follows immediately.
\begin{corollary}
For any pair $(a,b)$, there is at most one pair $(c,d)$ such that there is perfect state transfer from $(a,b)$ to $(c,d)$.\qed
\end{corollary}

Being periodic is a necessary condition for a state to have perfect state transfer and the period of a state involved in perfect state transfer can tell us the exact time when the perfect state transfer occurs. So periodicity of states provides a useful tool for analysis of perfect state transfer.
\begin{theorem}
[the Ratio Condition] Let $U(t)$ be the transition matrix corresponding to a graph $G$. Let $ \sum_r \theta_rE_r$ be the spectral decomposition of the Laplacian of $G$. Then $e_a-e_b$ is periodic in $G$ if and only if \[
\cfrac{\theta_r-\theta_s}{\theta_k-\theta_l}\in \mathbb{Q}\] for any $\theta_r,\theta_s,\theta_l,\theta_s$ in the eigenvalue support of $(e_a-e_b)$ with $\theta_l\neq \theta_k$.\qed
\end{theorem}

The following theorem can be viewed as a corollary of the ratio condition. The original proof can be found in Coutinho and Godsil~\cite{quantumnote} stated in terms of vertex states. A detailed proof using pair state transfer can be found in \cite{CQTmmath}.
\begin{theorem}
\label{pedev}
Let $G$ be a graph with the Laplacian matrix $L$ and let $(a,b)$ be a pair of vertices of $G$ with eigenvalue support $S$. Then $e_a-e_b$ is periodic in $G$ if and only if either:
\begin{enumerate}[label=(\roman*)]
\item All the eigenvalues in $S$ are integers;
\item There is a square-free integer $\Delta$ such that all eigenvalues in $S$ are quadratic integers in $\mathbb{Q}(\sqrt{\Delta})$, and the difference of any two eigenvalues in $S$ is an integer multiple of $\sqrt{\Delta}$.\qed
\end{enumerate}
\end{theorem}

\begin{corollary}
\label{pedeigendiff}
If $e_a-e_b$ is periodic in graph $G$, then any two distinct eigenvalues in the eigenvalue supports of $e_a-e_b$ differ by at least one.\qed
\end{corollary}

\begin{corollary}
If a pair state is periodic in graph $G$ with period $\tau$, then $\tau\leq 2\pi$.\qed
\end{corollary}
Using the ratio condition, we give the following necessary and sufficient conditions for perfect pair state transfer. Due to the similarity of the proof, we omit the proof here.  One can refer to \cite[Theorem~2.4.4]{Coutinho2014a} for details.
\begin{theorem}
Let $X$ be a graph. Then $X$ admits perfect pair state transfer between $e_a-e_b$ and $e_c-e_d$ if and only if all of the following conditions hold.\begin{enumerate}[label=(\roman*)]
\item Pair states $e_a-e_b$ and $e_c-e_d$ are strongly cospectral. Let $\theta_0$ be an eigenvalue in $\Lambda_{ab,cd}^+$.
\item The eigenvalues in the eigenvalue support of $e_a-e_b$ are either all integers or all quadratic integers. Moreover, there is a square-free integer $\Delta$ such that all eigenvalues in the eigenvalue support are quadratic integers in $\mathbb{Q}(\sqrt{\Delta})$, and the difference of any two eigenvalues in the eigenvalue support is an integer multiple of $\sqrt{\Delta}$.
\item Let $g=\gcd\left( \{\frac{\theta_0-\theta_r}{\sqrt{\Delta}}\}_{r=0}^k\right)$. Then \begin{enumerate}[label=(\alph*)]
\item $\theta_r\in\Lambda_{ab,cd}^+$ if and only if $\frac{\theta_0-\theta_r}{g\sqrt{\Delta}}$ is even, and 
\item $\theta_r\in\Lambda_{ab,cd}^-$ if and only if $\frac{\theta_0-\theta_r}{g\sqrt{\Delta}}$ is odd.\qed
\end{enumerate}
\end{enumerate}
\end{theorem}

\subsection{Fixed Pair States}
Let $v$ be a vertex of graph $G$. We use $N(v)$ to denote the neighbours of $v$ in $G$. We say a pair state $e_a-e_b$ is fixed if for all non-negative $t\in \mathbb{R}$, \[U(t)(e_a-e_b)=\gamma(e_a-e_b)\] with $\gamma$ being a norm-one complex scalar. We prove that a pair $(a,b)$ is fixed  if and only if vertices $a,b$ are twins in $G$, which means that \[N(a)\backslash\{b\}=N(b)\backslash\{a\}.\] That a pair state $e_a-e_b$ is fixed implies that it can never have perfect pair state transfer. Notice that a fixed state can be viewed as a state that is periodic at $t$ for any non-negative $t\in \mathbb{R}$.

\begin{lemma}
The pair state $e_a-e_b$ is fixed in $G$ if and only if the density matrix of $e_a-e_b$ and the Laplacian matrix $L$ of $G$ commute.
\end{lemma}

\proof
Let $D$ denote the density matrix of $e_a-e_b$. For any non-negative $t$, we have \[
U(t)DU(-t) = D\] if and only if \[U(t)D=DU(t),\] which means that $D$ commutes with $U(t)$ for all $t$. This is equivalent to that $D$ commutes with $L$ as $U(t)=\sum_m^\infty \frac{(it)^m}{m!}L$.\qed

\begin{lemma}
Let $D$ denote the density matrix of $e_a-e_b$ and let $L$ denote the Laplacian matrix of graph $G$. Then $LD=DL$ if and only if vertices $a,b$ are twins in $G$.
\end{lemma}

\proof 
We know that $DL=LD$ if and only $LD$ is symmetric, which is equivalent to $N(a)=N(b)$.\qed

The theorem below follows immediately.

\begin{theorem}
The pair state $e_a-e_b$ is fixed if and only if vertices $a,b$ are twins in $G$.\qed
\end{theorem}

\begin{corollary}
Let $a,b$ be two vertices in $G$. If $N(a)=N(b)$, then $e_a-e_b$ has no perfect pair state transfer.
\end{corollary}

\proof 
From previous theorem, we know that $e_a-e_b$ is fixed.\qed

\begin{theorem}
Pair states $e_a-e_b$ are fixed if and only if there is only one eigenvalue in the eigenvalue support of $e_a-e_b$ and the eigenvalue must be an integer.
\end{theorem}

\proof 
Let $\sum_r \theta_rE_r$ denote the spectral decomposition of the Laplacian of the graph $G$. For any time $t$, we have
\[
U(t)(e_a-e_b)=\sum_r \theta_r E_r(e_a-e_b)=\gamma(e_a-e_b)
\]  
for some complex scalar $\gamma$ with $\abs{\gamma}=1$. This is equivalent to the assumption that for all eigenvalues $\theta_r$, we have 
\[
e^{it\theta_r}E_r(e_a-e_b)= \gamma E_r(e_a-e_b),
\] 
which gives us that
\[
\gamma= e^{it\theta_r}
\] 
for all $\theta_r$ at any time $t$. 

It follows that $e_a-e_b$ is fixed if and only if all the eigenvalues in the eigenvalue support coincide. If $\theta_r$ is an eigenvalue in the eigenvalue support, then all the algebraic conjugates of $\theta_r$ are also in the eigenvalue suppport. So we can conclude that $e_a-e_b$ is fixed if and only if the eigenvalue support of $e_a-e_b$ is $\{ \theta\}$ for some integer $\theta$.\qed

\begin{corollary}
In a graph $G$, vertices $a$ and $b$ are twins if and only if the eigenvalue support of $e_a-e_b$ consists of one integer eigenvalue.\qed
\end{corollary} 

This is a feature that distinguishes vertex state transfer and pair state transfer.  In a 
connected graph with at least two vertices, the eigenvalue support of a vertex state must have 
size at least two, while the eigenvalue support of a pair state can have size one.

\section{Algebraic Properties}
\label{algebraic properties}

In this section, we show how algebraic properties of the underlying graphs can help us to get more information about state transfer.

\begin{theorem}
\label{aut same eigensupp} Let $G$ be a graph.
If there is a permutation $\sigma\in \aut{G}$ such that $\sigma(e_a-e_b)=e_c-e_d$, then $e_a-e_b$ and $e_c-e_d$ have the same eigenvalue support.
\end{theorem}

\proof 
Let $P$ denote the permutation matrix associated with $\sigma\in \aut{G}$. Since $\col(A)$ is invariant under $P$, we have that \[PA=AP.\] Let $\Delta$ denote the degree matrix of $G$. We know that $\Delta$ is a diagonal matrix, so that $P$ commutes with $\Delta$. Let $L$ denote the Laplacian matrix of $G$. Thus, we have that \[ LP=\left(\Delta-A\right)P=P\left(\Delta-A\right)=PL.\] Let $\sum_r\theta_r E_r$ denote the spectral decomposition of $L$.  Since $E_r$ is a polynomial in $L$ and hence, we know that $L$ commutes with $E_r$. Then we have that \[
PE_r(e_a-e_b)=E_rP(e_a-e_b)=E_r(e_c-e_d).\]
We know that $\theta_r$ is not in the eigenvalue support of $e_a-e_b$ if and only if \[E_r(e_a-e_b) = 0.\] Since $P$ acts on $E_r(e_a-e_b)$ by permuting its entries, we can see that $E_r(e_a-e_b) = 0$ if and only if \[
PE_r(e_a-e_b)=E_r(e_c-e_d)=0.\] Thus, we can conclude that $\theta_r$ is not in the eigenvalue support of $e_a-e_b$ if and only if $\theta_r$ is not in the eigenvalue support of $e_c-e_d$.\qed
One immediate consequence is that all the edge states in an edge-transitive graph have the same eigenvalue support. Actually the eigenvalue support of an edge state of a edge-transitive graph is consist of all the non-zero eigenvalues.

\begin{theorem}
\label{all nonzero eigenvalue in the eigensupp}
If $G$ is a connected edge-transitive graph, then the eigenvalue support of an edge state of $G$ consists of all the non-zero eigenvalues.
\end{theorem}

\proof Let $\sum_r \theta_r E_r$ denote the spectral decomposition of the Laplacian matrix of $G$. Assume towards contradictions that $E_s$ is the spectral idempotent corresponding to a non-zero eigenvalue $\theta_s$ such that \[
E_s(e_a-e_b)=0\] for all $(a,b)\in E(G)$. Since $G$ is connected, for any $u,v\in V(G)$, there exist a path $P$ from $u$ to $v$. Then for any two vertices $i,j$ on $P$, we must have \[E_s e_i= E_s e_j.\] We can conclude that all the columns of $E_s$ are equal, which contradict to that $\theta_s$ is a non-zero eigenvalue. Thus, we know that a non-zero eigenvalue $\theta_r$ must in the eigenvalue support of some edge state of $G$. Since $G$ is edge-transitive, we know that the eigenvalue supports of $e_a-e_b$ for all $(a,b)\in E(G)$ are equal. Therefore, all the non-zero eigenvalues are in eigenvalue support of edge states of $G$.\qed
In the proof of Theorem~\ref{aut same eigensupp}, we prove that \[LP=PL.\]
The transition matrix associated with graph $G$ is a polynomial in $L$. So for any permutation matrix $P$  from $\aut{G}$, we have \[
U(t)P=PU(t).\]
Using this, the following Lemma is proved in \cite[Corollary 9.2]{Godsil2012}.

\begin{lemma}
\label{same stabilizer}
If a graph $G$ admits perfect state transfer between $e_a-e_b$ to $e_c-e_d$, then the stabilizer of $e_a-e_b$ is the same as the stabilizer of $e_c-e_d$ in $\aut{G}$.\qed
\end{lemma}

\begin{lemma}
If graph $G$ admits perfect state transfer between $e_a-e_b$ to $e_c-e_d$, then all the pair states in the orbit of $e_a-e_b$ under $\aut{G}$ have perfect state transfer.
\end{lemma}

\proof Assume there exist time $\tau$ such that $U(\tau)(e_a-e_b)=\gamma(e_c-e_d)$ for some $\abs{\gamma}=1$. Let $P$ denote a permutation matrix associated with a $\sigma\in \aut{G}$ and $P(e_a-e_b)=e_{a'}-e_{b'}$. By Lemma~\ref{same stabilizer}, we know $P$ does not fix $(e_c-e_d)$ and assume $P(e_c-e_d)=e_{c'}-e_{d'}$. Then we have \begin{gather*}
U(\tau)P(e_a-e_b)=\gamma P(e_c-e_d)\\
U(\tau)(e_{a'}-e_{b'})=\gamma( e_{c'}-e_{d'}).
\end{gather*} Thus, there is also perfect state transfer between $e_{a'}-e_{b'}$ and $e_{c'}-e_{d'}$ at time $\tau$.\qed 

By the monogamy of perfect state transfer, the following results are immediate.

\begin{corollary}
If there is perfect state transfer between $e_a-e_b$ and $e_c-e_d$ in graph $G$, then the orbit of $e_a-e_b$ and the orbit of $e_c-e_d$ under $\aut{G}$ must have the same size.\qed
\end{corollary}

\begin{corollary}
Given an edge-transitive graph $G$, if perfect edge state transfer occurs in $G$, then all the edges have perfect state transfer.\qed
\end{corollary}

By monogamy property of perfect state transfer, we know that perfect state edge transfer in an edge-transitive graph partition edges into pairs. 
\begin{corollary}
\label{noPST in odd edge-transitive}
Let $G$ be an edge-transitive graph with $n$ edges. If $n$ is odd, there is no edge perfect state transfer in $G$.\qed
\end{corollary}

\section{Constructions}
In this section, we show how to use complements and Cartesian products to build infinite families of graphs with perfect pair state transfer.

\subsection{Complements}
\label{complement}
We use standard algebraic graph theory result to show that complementation preserves perfect pair state transfer. Let $\comp{G}$ denote the complement of a graph $G$.

\begin{lemma}
\label{disconnected Leigenvec}
Let $G$ be a graph with $n$ vertices and $L$ denote the Laplacian matrix of $G$. Then every Laplacian eigenvector of $G$ with non-zero eigenvalue $\theta$ is a Laplacian eigenvector of $\comp{G}$ with eigenvalue $n-\theta$.\qed
\end{lemma}

\begin{theorem}
\label{pstcomple}
There is perfect state transfer between $(e_a-e_b)$ and $(e_c-e_d)$ in graph $G$ if and only if there is perfect state transfer between $(e_a-e_b)$ and $(e_c-e_d)$ in $\comp{G}$.
\end{theorem}

\proof
Let $S=\{\theta_1,\theta_2,\cdots,\theta_r\}$ denote the eigenvalue support of $(e_a-e_b)$ 
and $(e_c-e_d)$ in $G$ and $\sum_r \theta_r E_r$ denote the spectral decomposition of the 
Laplacian of $G$. Let \[a_j = (E_j)_{ac} +(E_j)_{ad}- (E_j)_{bc}+(E_j)_{bd}\] for all 
eigenvalues $\theta_j\in S$. Then we have that
\begin{align*}
&\abs{\frac{1}{2}(e_c-e_d)^T U(t)(e_a-e_b)}^2\\
&=\abs{\frac{1}{2} \sum_{j=1}^r e^{it\theta_j}\left((E_j)_{ac} +(E_j)_{ad}
	-(E_j)_{bc}+(E_j)_{bd}\right)}^2\\
&=\abs{\frac{1}{2}\left(a_1e^{it\theta_1}+a_2e^{it\theta_2}+\cdots+a_re^{it\theta_r}\right)}^2\\
&=\frac{1}{4}\big( \left(a_1\cos(\theta_1 t)+a_2\cos(\theta_2 t)+\cdots+a_r\cos(\theta_r t)\right)^2\\
&\phantom{=}+\left(a_1\sin(\theta_1 t)+a_2\sin(\theta_2 t)+\cdots+a_r\sin(\theta_r t)\right)^2\big)\\
&=\frac{1}{4}\left( a_1^2+a_2^2+\cdots+a_r^2+\sum_{r\neq s} 2a_ra_s\cos((\theta_r-\theta_s)t)\right)
\end{align*}
By Lemma~\ref{disconnected Leigenvec}, we know that the eigenvalue support $\comp{S}$ of $(e_a-e_b)$ and $(e_c-e_d)$  in $\comp{G}$ is $\{n-\theta_1,n-\theta_2,\cdots,n-\theta_r\}$. Since zero is never in the eigenvalue support, the spectral idempotent $\comp{E_r}$ of the Laplacian $\comp{L}$ of $\comp{G}$ with eigenvalue $n-\theta_r$ is the same as $E_r$ with eigenvalue $\theta_r$ of $L$ for all eigenvalues in the eigenvalue support of $e_a-e_b$ in $G$.
Let $\comp{U}(t) = \exp(it\comp{L})$ be the transition matrix associated with $\comp{G}$. We have that 
\begin{align*}
&\abs{\frac{1}{2}(e_c-e_d)^T \comp{U}(t) (e_a-e_b)}^2\\
&=\abs{\frac{1}{2} \sum_{j=1}^r e^{it(n-\theta_j)}\left((E_j)_{ac} +(E_j)_{ad}- (E_j)_{bc}+(E_j)_{bd}\right)}^2 \\
&=\frac{1}{4}\left( a_1^2+a_2^2+\cdots+a_r^2+\sum_{r\neq s} 2a_ra_s\cos\left( (n-\theta_r)t-(n-\theta_s)t \right)\right)\\
&=\frac{1}{4}\left( a_1^2+a_2^2+\cdots+a_r^2+\sum_{r\neq s} 2a_ra_s\cos\left( (\theta_s-\theta_r)t \right)\right)
\end{align*}
Since cosine is an even function, we get that \[
\abs{\frac{1}{2}(e_c-e_d)^T \comp{U}(t) (e_a-e_b)}^2 = 
\abs{\frac{1}{2}(e_c-e_d)^T U(t) (e_a-e_b)}^2.\] Therefore, there is perfect state transfer between $(e_a-e_b)$ and $(e_c-e_d)$ in graph $G$ if and only if there is perfect state transfer between them in the complement of $G$.\qed
Let $G_1,G_2$ be two graphs. Let $E'$ denote the set of all the edges with one end in $V(G_1)$ and the other end in $V(G_2)$. The join graph of $G_1$ and $G_2$ is a graph $G$ such that \[
V(G)=V(G_1)\cup V(G_2), E(G)= E(G_1)\cup E(G_2)\cup E' .\]
\begin{corollary} Let $G$ be a graph and $a,b,c,d$ are vertices in $G$. There is perfect state transfer between $e_a-e_b$ and $e_c-e_d$ in $G$ if and only if there is perfect state transfer between $e_a-e_b$ and $e_c-e_d$ in the join graph of $G$ and $H$ for a graph $H$.
\end{corollary}
Notice that when $H$ is a simple graph with one vertex, the join graph of a graph $G$ and $H$ is a cone graph $G$. So we can see that if there is perfect pair state transfer in a graph $G$, using Theorem~\ref{pstcomple}, we can easily construct a cone graph of $G$ to obtain a new graph that admits perfect pair state transfer.

Theorem~\ref{pstcomple} also allows us to characterize perfect state transfer in some graphs with special structure.
\begin{corollary}
\label{Kn-1 pst}
Let $K_n$ be a complete graph on $n$ vertices and $V(K_n)=\{v_1,v_2,\cdots,v_n\}$. Let $G$ denote the graph obtained from $K_n$ by deleting edge $(v_1,v_2)$. Then there is perfect state transfer between $e_1-e_i$ and $e_2-e_i$ for all $i\in\{3,4,\cdots,n\}$.
\end{corollary}
\subsection{Cartesian Products}
Cartesian product is also an operator that we can use to construct infinite families of graphs with perfect pair state transfer.

Let $G,H$ be two graphs, their Cartesian product has vertex set $V(G)\times V(H)$, where $(g_1,h_1)$ is adjacent to $(g_2,h_2)$ if and only if either
\begin{enumerate}[label=(\roman*)]
\item $g_1= g_2$ in $G$ and $h_1$ is adjacent to $h_2$ in $H$, or
\item $g_1$ is adjacent to $g_2$ in $G$ and $h_1=h_2$ in $H$.
\end{enumerate}

\begin{lemma}
Let $G,H$ be graphs with Laplacian matrices $L_G$ of order $n\times n$, $L_H$ of order $m\times m$ respectively. Let $G\square H$ denote the Cartesian product of $G$ and $H$ with the Laplacian matrix $L_{G\square H}$.  Then $L_{G\square H} = L_G\otimes I + I\otimes L_H.$\qed
\end{lemma}

\begin{lemma}
Let $G,H$ be two graphs with transition matrices $U_G(t) = \exp(itL_G)$ and $U_H(t) = \exp(itL_H)$ respectively. Let $U_{G\square H} (t)= \exp(itL_{G\square H})$ denote the transition matrix of $G\square H$. Then $U_{G\square H}(t) = U_G(t) \otimes U_H(t)$.
\end{lemma}
\proof Let $L_G$ be a matrix of order $n\times n$ and let $L_H$ be a matrix of order $m\times m$. If $M$ is a matrix of order $m$ and $N$ is a matrix of order $n\times n$, the Kronecker sum of $M$ and $N$ is \[M\oplus N = M\otimes I_n+I_m\otimes N.\]
Using the Kronecker sum and previous lemma, we have
\begin{align*}
U_{G\square H}(t) & = \exp(itL_{G\square H})\\
&=\exp\left(it( L_G\otimes I_m + I_n\otimes L_H)\right)\\
&=\exp\left( it (L_G \oplus L_H)\right)\\
&=\exp(itL_G)\otimes\exp(itL_H)\\
&=U_G(t)\otimes U_H(t).
\end{align*}\qed

\begin{theorem}
\label{catpairpst}
Let $G,H$ be two graphs, let $(a,b),(c,d)$ be two pairs of vertices in $G$ and let $(\alpha,\beta),(\gamma,\kappa)$ be two pairs of vertices in $H$. There is perfect state transfer between the pair $\{(a,\alpha),(b,\beta)\}$  and the pair $\{(c,\gamma),(d,\kappa)\}$ in $G\square H$ at time $t$ if and only if both of the following conditions hold:
\begin{enumerate}[label=(\roman*)]
\item There is perfect pair state transfer between the pair $(a,b)$ and $
(c,d)$ in $G$ at time $t$.
\item There is perfect pair state transfer between edges $(\alpha,\beta)$ and $(\gamma,\kappa)$ in $H$ at time $t$.
\end{enumerate}
\end{theorem}

\proof
The state associated with the pair $\{(a,\alpha),(b,\beta)\}$ is \[
\cfrac{1}{2}\left(( e_a-e_b)\otimes (e_\alpha -e_\beta)\right)\] and then we can see that the density matrix of this edge is \[D_{ab}\otimes D_{\alpha\beta}.\] Similarly, the density matrix of the edge $\{(c,\gamma),(d,\kappa)\}$ is $D_{cd}\otimes D_{\gamma\kappa}$. There is perfect state transfer between $\{(a,\alpha),(b,\beta)\}$  and $\{(c,\gamma),(d,\kappa)\}$ at time $t$ if and only if  
\[ 
	U_{G\square H}(t)\cdot D_{ab}\otimes D_{\alpha\beta}\cdot U_{G\square H}(-t) 
		= D_{cd}\otimes D_{\gamma\kappa}.
\] 
By the previous corollary, we have that
\begin{align*}
U_{G\square H}(t)&\cdot D_{ab}\otimes D_{\alpha\beta} \cdot U_{G\square H}(-t)\\
 &= U_{G\square H}(t)\cdot D_{ab}\otimes D_{\alpha\beta}\cdot\left( U_{G}(-t) \otimes U_H(-t) \right) \\
 &=\left( U_{G}(t) \otimes U_H(t) \right)\cdot(D_{ab}U_G(-t)\otimes D_{\alpha\beta} U_H(-t)\\
 &=\left(U_G(t)D_{ab}U_G(-t)\right)\otimes \left(U_H(t)D_{\alpha\beta} U_H(-t)\right)\\
 &=D_{cd}\otimes D_{\gamma\kappa},
\end{align*} which is equivalent to that there is perfect Laplacian state transfer between $(a,b)$,$(c,d)$ in $G$ at time $t$ and at the same time there is perfect state transfer between edge $(\alpha,\beta)$ and $(\gamma,\kappa)$ in $H$ .\qed

Notice that in the case that the pair $\{(a,\alpha),(b,\beta)\}$ and the pair $\{(c,\gamma),(d,\kappa)\}$ are both edges, which means that $a=b$ and $c=d$, perfect pair state transfer in $G\square H$ is a combination of perfect vertex state transfer and perfect pair state transfer.
 
\begin{corollary}
\label{catpst}
Let $G,H$ be two graphs, let $a,b$ be two vertices in $G$ and let $\alpha,\beta,\gamma,\kappa$ be vertices in $H$. There is perfect state transfer between the edge $\{(a,\alpha),(a,\beta)\}$  and the edge $\{(b,\gamma),(b,\kappa)\}$ in $G\square H$ at time $t$ if and only if both of the following conditions hold:
\begin{enumerate}[label=(\roman*)]
\item There is perfect Laplacian vertex state transfer between vertices $a$ and $b$ in $G$ at time $t$.
\item There is perfect pair state transfer between edges $(\alpha,\beta)$ and $(\gamma,\kappa)$ in $H$ at time $t$.
\end{enumerate}
\end{corollary}

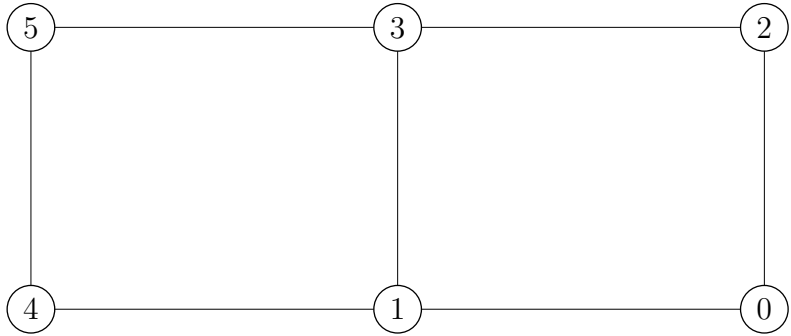
\begin{figure}[h]
\begin{center}
\begin{tikzpicture}[scale=0.75]

\definecolor{cv0}{rgb}{0.0,0.0,0.0}
\definecolor{cfv0}{rgb}{1.0,1.0,1.0}
\definecolor{clv0}{rgb}{0.0,0.0,0.0}
\definecolor{cv1}{rgb}{0.0,0.0,0.0}
\definecolor{cfv1}{rgb}{1.0,1.0,1.0}
\definecolor{clv1}{rgb}{0.0,0.0,0.0}
\definecolor{cv2}{rgb}{0.0,0.0,0.0}
\definecolor{cfv2}{rgb}{1.0,1.0,1.0}
\definecolor{clv2}{rgb}{0.0,0.0,0.0}
\definecolor{cv3}{rgb}{0.0,0.0,0.0}
\definecolor{cfv3}{rgb}{1.0,1.0,1.0}
\definecolor{clv3}{rgb}{0.0,0.0,0.0}
\definecolor{cv4}{rgb}{0.0,0.0,0.0}
\definecolor{cfv4}{rgb}{1.0,1.0,1.0}
\definecolor{clv4}{rgb}{0.0,0.0,0.0}
\definecolor{cv5}{rgb}{0.0,0.0,0.0}
\definecolor{cfv5}{rgb}{1.0,1.0,1.0}
\definecolor{clv5}{rgb}{0.0,0.0,0.0}
\definecolor{cv0v1}{rgb}{0.0,0.0,0.0}
\definecolor{cv0v2}{rgb}{0.0,0.0,0.0}
\definecolor{cv1v3}{rgb}{0.0,0.0,0.0}
\definecolor{cv1v4}{rgb}{0.0,0.0,0.0}
\definecolor{cv2v3}{rgb}{0.0,0.0,0.0}
\definecolor{cv3v5}{rgb}{0.0,0.0,0.0}
\definecolor{cv4v5}{rgb}{0.0,0.0,0.0}
\Vertex[style={minimum size=1.0cm,draw=cv0,fill=cfv0,text=clv0,shape=circle},LabelOut=false,L=\hbox{$0$},x=13cm,y=0cm]{v0}
\Vertex[style={minimum size=1.0cm,draw=cv1,fill=cfv1,text=clv1,shape=circle},LabelOut=false,L=\hbox{$1$},x=6.5cm,y=0cm]{v1}
\Vertex[style={minimum size=1.0cm,draw=cv2,fill=cfv2,text=clv2,shape=circle},LabelOut=false,L=\hbox{$2$},x=13cm,y=5.0cm]{v2}
\Vertex[style={minimum size=1.0cm,draw=cv3,fill=cfv3,text=clv3,shape=circle},LabelOut=false,L=\hbox{$3$},x=6.5cm,y=5cm]{v3}
\Vertex[style={minimum size=1.0cm,draw=cv4,fill=cfv4,text=clv4,shape=circle},LabelOut=false,L=\hbox{$4$},x=0cm,y=0cm]{v4}
\Vertex[style={minimum size=1.0cm,draw=cv5,fill=cfv5,text=clv5,shape=circle},LabelOut=false,L=\hbox{$5$},x=0cm,y=5cm]{v5}
\Edge[lw=0.01cm,style={color=cv0v1,},](v0)(v1)
\Edge[lw=0.01cm,style={color=cv0v2,},](v0)(v2)
\Edge[lw=0.01cm,style={color=cv1v3,},](v1)(v3)
\Edge[lw=0.01cm,style={color=cv1v4,},](v1)(v4)
\Edge[lw=0.01cm,style={color=cv2v3,},](v2)(v3)
\Edge[lw=0.01cm,style={color=cv3v5,},](v3)(v5)
\Edge[lw=0.01cm,style={color=cv4v5,},](v4)(v5)
\end{tikzpicture}
\caption{$P_2\square P_3$}
\label{p23}
\end{center}
\end{figure}

The example given by Coutinho in \cite[Section~2.4]{Coutinho2014a} shows that $P_2$ admits perfect state transfer with respect to its Laplacian matrix between its vertices at time $\frac{\pi}{2}$. As the example shown in Section~\ref{trans mtx}, there is perfect pair state transfer between its edges in $P_3$ at time $\frac{\pi}{2}$. By Theorem~\ref{catpst}, there is perfect state transfer from $e_3-e_5$ to $e_1-e_0$ and from $e_2-e_3$ to $e_1-e_4$ in Figure~\ref{p23}. 
\begin{figure}[H]
\begin{center}
\begin{tikzpicture}[scale=0.75]
\definecolor{cv0}{rgb}{0.0,0.0,0.0}
\definecolor{cfv0}{rgb}{1.0,1.0,1.0}
\definecolor{clv0}{rgb}{0.0,0.0,0.0}
\definecolor{cv1}{rgb}{0.0,0.0,0.0}
\definecolor{cfv1}{rgb}{1.0,1.0,1.0}
\definecolor{clv1}{rgb}{0.0,0.0,0.0}
\definecolor{cv2}{rgb}{0.0,0.0,0.0}
\definecolor{cfv2}{rgb}{1.0,1.0,1.0}
\definecolor{clv2}{rgb}{0.0,0.0,0.0}
\definecolor{cv3}{rgb}{0.0,0.0,0.0}
\definecolor{cfv3}{rgb}{1.0,1.0,1.0}
\definecolor{clv3}{rgb}{0.0,0.0,0.0}
\definecolor{cv4}{rgb}{0.0,0.0,0.0}
\definecolor{cfv4}{rgb}{1.0,1.0,1.0}
\definecolor{clv4}{rgb}{0.0,0.0,0.0}
\definecolor{cv5}{rgb}{0.0,0.0,0.0}
\definecolor{cfv5}{rgb}{1.0,1.0,1.0}
\definecolor{clv5}{rgb}{0.0,0.0,0.0}
\definecolor{cv0v3}{rgb}{0.0,0.0,0.0}
\definecolor{cv0v4}{rgb}{0.0,0.0,0.0}
\definecolor{cv0v5}{rgb}{0.0,0.0,0.0}
\definecolor{cv1v2}{rgb}{0.0,0.0,0.0}
\definecolor{cv1v5}{rgb}{0.0,0.0,0.0}
\definecolor{cv2v4}{rgb}{0.0,0.0,0.0}
\definecolor{cv2v5}{rgb}{0.0,0.0,0.0}
\definecolor{cv3v4}{rgb}{0.0,0.0,0.0}
\Vertex[style={minimum size=1.0cm,draw=cv0,fill=cfv0,text=clv0,shape=circle},LabelOut=false,L=\hbox{$0$},x=7cm,y=5cm]{v0}
\Vertex[style={minimum size=1.0cm,draw=cv1,fill=cfv1,text=clv1,shape=circle},LabelOut=false,L=\hbox{$1$},x=-3cm,y=2.50cm]{v1}
\Vertex[style={minimum size=1.0cm,draw=cv2,fill=cfv2,text=clv2,shape=circle},LabelOut=false,L=\hbox{$2$},x=0cm,y=0cm]{v2}
\Vertex[style={minimum size=1.0cm,draw=cv3,fill=cfv3,text=clv3,shape=circle},LabelOut=false,L=\hbox{$3$},x=10cm,y=2.50cm]{v3}
\Vertex[style={minimum size=1.0cm,draw=cv4,fill=cfv4,text=clv4,shape=circle},LabelOut=false,L=\hbox{$4$},x=7.0cm,y=0cm]{v4}
\Vertex[style={minimum size=1.0cm,draw=cv5,fill=cfv5,text=clv5,shape=circle},LabelOut=false,L=\hbox{$5$},x=0.0cm,y=5cm]{v5}
\Edge[lw=0.01cm,style={color=cv0v3,},](v0)(v3)
\Edge[lw=0.01cm,style={color=cv0v4,},](v0)(v4)
\Edge[lw=0.01cm,style={color=cv0v5,},](v0)(v5)
\Edge[lw=0.01cm,style={color=cv1v2,},](v1)(v2)
\Edge[lw=0.01cm,style={color=cv1v5,},](v1)(v5)
\Edge[lw=0.01cm,style={color=cv2v4,},](v2)(v4)
\Edge[lw=0.01cm,style={color=cv2v5,},](v2)(v5)
\Edge[lw=0.01cm,style={color=cv3v4,},](v3)(v4)
\end{tikzpicture}
\caption{the complement of the graph in Figure~\ref{p23}}
\label{compp23}
\end{center}
\end{figure}
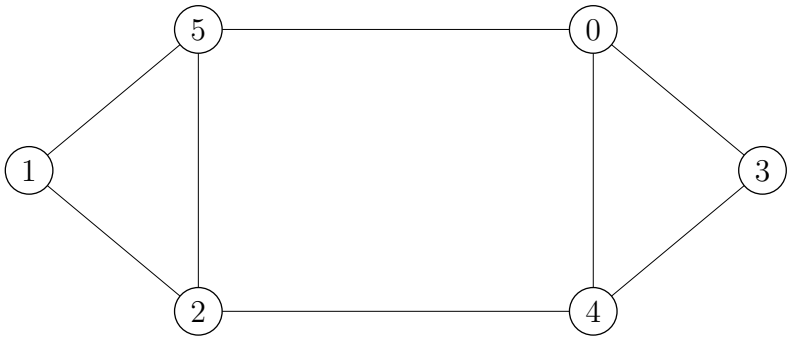
\section{Transitivity}
\label{transitivity}

The graph in Figure~\ref{compp23} is the complement of the graph in Figure~\ref{p23}. We know that there is perfect pair state transfer from $e_3-e_5$ to $e_1-e_0$ and also from $e_2-e_3$ to $e_1-e_4$ in the graph shown in Figure~\ref{p23}. By Theorem~\ref{pstcomple}, we know that there are also perfect pair states transfer from $e_3-e_5$ to $e_1-e_0$ and from $e_2-e_3$ to $e_1-e_4$ the graph shown in Figure~\ref{compp23}. 

Actually the graph in Figure~\ref{compp23} also admits perfect pair state transfer between $e_2-e_5$ and $e_0-e_4$. This is what we call ``the transitivity phenomenon". This phenomenon can never happen in the case of vertex state transfer due to the monogamy and symmetry properties of perfect state transfer.

\begin{theorem}
\label{transitivitythm}
Suppose there is perfect state transfer between $e_a-e_b$ and $e_\alpha-e_\beta$ at time $\tau$ 
in $G$ and there is also perfect state transfer between $e_b-e_c$ and $e_\beta-e_\gamma$ at the 
same time $\tau$ in $G$. Then there is perfect state transfer between $e_a-e_c$ and 
$e_\alpha-e_\gamma$ at time $\tau$ in $G$.
\end{theorem}

\proof
Let $D_{ab}$ denote the density matrix of $e_a-e_b$ and $D_{bc}$ denote the density matrix of $e_b-e_c$. We have that \[
D_{ab}=\frac{1}{2}(e_a-e_b)(e_a-e_b)^T,\quad D_{bc}=\frac{1}{2}(e_b-e_c)(e_b-e_c)^T.\] Using that \[(e_a-e_b)^T(e_b-e_c) = (e_b-e_c)^T(e_a-e_b) =-1,\]
we can write the density matrix of $e_a-e_c$ in terms of $D_{ab}$ and $D_{bc}$ in the following way:
\begin{align*}
D_{ac} &= \frac{1}{2}(e_a-e_c)(e_a-e_c)^T\\
&=\frac{1}{2}\left( (e_a-e_b)+(e_b-e_c)\right)\left( (e_a-e_b)+(e_b-e_c)\right)^T\\
&=\frac{1}{2}\big((e_a-e_b)(e_a-e_b)^T+(e_b-e_c)(e_b-e_c)^T \\&
\phantom{=}+ (e_a-e_b)(e_b-e_c)^T+(e_b-e_c)(e_a-e_b)^T\big)\\
&=\frac{1}{2}\big((e_a-e_b)(e_a-e_b)^T+(e_b-e_c)(e_b-e_c)^T\\
&\phantom{=}-(e_a-e_b)(e_a-e_b)^T(e_b-e_c)(e_b-e_c)^T\\
&\phantom{=}-(e_b-e_c)(e_b-e_c)^T(e_a-e_b)(e_a-e_b)^T\big)\\
&= D_{ab}+D_{bc}-2D_{ab}D_{bc}-2D_{bc}D_{ab} 
\end{align*}
As the above shows, we have 
\[
D_{ac} = D_{ab}+D_{bc}-2D_{ab}D_{bc}-2D_{bc}D_{ab}.
\] 
Similarly, we have 
\[  
	D_{\alpha\gamma} = D_{\alpha\beta}+D_{\beta\gamma}-2D_{\alpha\beta}D_{\beta\gamma}
		-2D_{\beta\gamma}D_{\alpha\beta}.
\]

Now consider $U(\tau) D_{ac} U(-\tau)$. Since we know that 
\[
U(\tau) D_{ab} U(-\tau)=D_{\alpha\beta} \quad \text{and}\quad U(\tau) D_{bc} U(-\tau)=D_{\beta\gamma},
\] 
we have
\begin{align*}
U(\tau) D_{ac} U(-\tau) &= U(\tau)\left(D_{ab}+D_{bc}-2D_{ab}D_{bc}-2D_{bc}D_{ab}\right) U(-\tau) \\
&= D_{\alpha\beta}+D_{\beta\gamma} -2 U(\tau) D_{ab}D_{bc} U(-\tau)-2U(\tau) D_{bc}D_{ab} U(-\tau).
\end{align*}
Using $U(-\tau)\cdot U(\tau) = 1$, we get
\[
	U(\tau) D_{ab}D_{bc} U(-\tau) = U(\tau) D_{ab}U(-\tau)\cdot U(\tau)D_{bc} U(-\tau) 
		=D_{\alpha\beta}D_{\beta\gamma}
\] 
and similarly,
\[
U(\tau) D_{bc}D_{ab} U(-\tau) =U(\tau) D_{bc}U(-\tau)\cdot U(\tau)D_{ab} U(-\tau) =D_{\beta\gamma}D_{\alpha\beta}.
\] 
Thus, we get that 
\[
U(\tau) D_{ac} U(-\tau) = D_{\alpha\beta}+D_{\beta\gamma}-2D_{\alpha\beta}D_{\beta\gamma}-2D_{\beta\gamma}D_{\alpha\beta} = D_{\alpha\gamma}.
\]
Therefore, there is perfect state transfer between $e_a-e_c$ and $e_\alpha-e_\gamma$ at time $\tau$.\qed

\section{Special Classes}
This section we discuss pair state transfer on paths and cycles. Since we have proved that perfect pair state transfer are equivalent up to taking complements of underlying graphs, here we exclude the case of perfect pair state transfer between pairs that are both non-edges.
 
We show that $C_4$ is the only cycle and $P_3,P_4$ are the only paths that have perfect pair state transfer. We observe an interesting correspondence of perfect state transfer between graphs and their line graphs when graphs are paths and cycles.
\label{specialcases}

\subsection{Cycles}
We use $C_n$ to denote the cycle on $n$ vertices and $A(C_n), L(C_n)$ to denote the adjacency and Laplacian matrix of $C_n$ respectively.

Since here we only consider the case when at least one of the pairs that has perfect pair state transfer is an edge, we use a bound on $n$ such that $C_n$ can have a periodic edge state to eliminate the cases when $C_n$ can have perfect pair state transfer. We show that $C_4$ is the only cycle that has perfect pair state transfer.

\begin{lemma}
\label{cyceval}
Laplacian eigenvectors of $C_n$ are \[v_k=
\begin{pmatrix}
1\\ \omega^k\\\omega^{2k}\\ \omega^{3k}\\\vdots\\ \omega^{(n-1)k}\end{pmatrix}\] for $k=0,1,\cdots,n-1$ where $\omega=e^{\frac{2\pi}{n}i}$ with eigenvalues \[2-2\cos\frac{2\pi k}{n}\] for $k = 0,1,\cdots,n-1$.\qed
\end{lemma}

Since we have \[
\cos\frac{2\pi (n-r)}{n}=\cos\left(2\pi-\frac{2\pi r}{n}\right)=\cos\frac{2\pi r}{n}, \] we know that $k=r$ and $k=n-r$ produce the same eigenvalue for $r\in\{1,2,\cdots,n-1\}$. Thus, we can conclude that $L(C_n)$ has $\lfloor \frac{n}{2}\rfloor$ distinct non-zero eigenvalues. 

Using Theorem~\ref{all nonzero eigenvalue in the eigensupp}, the Lemma below follows immediately.
\begin{lemma}
\label{cyceigensuppsize}
Every edge state of $C_n$ has eigenvalue support of size $\lfloor \frac{n}{2}\rfloor$.\qed
\end{lemma}

\begin{theorem}
\label{cycle edge pst}
There is perfect pair state transfer in $C_n$ if and only if $n=4$. 
\end{theorem}
\proof By Lemma~\ref{cyceval}, we know that the Laplacian eigenvalues of $C_n$ are \[
0\leq 2-2\cos\left(\frac{2\pi k}{n}\right)\leq 4\] for $k=0,1,\cdots,n-1$. By Corollary~\ref{pedeigendiff}, we know that for  an edge state to be periodic, the size of eigenvalue support must be at most $4$. Then by Lemma~\ref{cyceigensuppsize}, we know that for $C_n$ to have a periodic edge state, we must have $
3\leq n\leq 9$.

Using Theorem~\ref{pedev}, we can find that there are no periodic edge states in $C_n$ when $n=7,8,9$ which implies that there is no perfect edge state transfer in $C_n$. Since cycles are edge-transitive, by Corollary~\ref{noPST in odd edge-transitive}, we know there is no perfect state transfer in $C_3$ and $C_5$. 

Computing \[
\abs{\frac{1}{2}(e_a - e_b)^TU(t)(e_c - e_d)} ^2\] for all vertex-pairs $(a,b),(c,d)$ in $V(C_n)$ when $n=4,6$, we can conclude that the only cycle that has perfect pair state transfer is $C_4$.\qed

At time $\frac{\pi}{2}$, there is perfect state transfer between the opposite edges in $C_4$.
\subsection{Paths}
Let $P_n$ denote the path on $n$ vertices such that $V(P_n)=\{1,2,\cdots,n\}$. We show that $P_3,P_4$ are the only two paths where perfect pair state transfer occurs.
\begin{lemma}
\label{pathevt}
The Laplacian eigenvector with eigenvalue $2-2\cos\frac{\pi r}{n}$ of $P_n$ is \[ 2\sin\left(\frac{r\pi}{2n}\right)
\begin{pmatrix}
\cos\left(1\frac{r\pi}{2n}\right)\\
\cos\left(3\frac{r\pi}{2n}\right)\\
\cos\left(5\frac{r\pi}{2n}\right)\\
\vdots\\
\cos\left((2n-3)\frac{r\pi}{2n}\right)\\
\cos\left((2n-1)\frac{r\pi}{2n}\right)
\end{pmatrix}\] for $r = 0,1,\cdots,n-1$.\qed
\end{lemma} 
Using automorphisms of path graphs and Theorem~\ref{aut same eigensupp}, we can prove the symmetry of the eigenvalue supports of the edge states of $P_n$.
\begin{lemma}
\label{pathsymmetric eigen_supp}
Let $(k,k+1)$ be an edge of $P_n$ with $1\leq k\leq n-1$. Then the eigenvalue supports of the edge states associated with $(k,k+1)$ and $(n-k,n-k+1)$ are the same.\qed
\end{lemma}

\begin{lemma}
\label{path eigen_supp size}
Let $S$ denote the eigenvalue support of an edge state in $P_n$. Then \[ \abs{S} \geq \frac{n}{2}.\]
\end{lemma}
\proof We want to prove that there are at most ${n}/{2}$ eigenvalues that are not in the eigenvalue support of an edge state in $P_n$. 

Let $E_r$ denote the spectral idempotent of $L(P_n)$ with eigenvalue $2-2\cos\left(\pi r/{n}\right)$. Since $0$ is never in the eigenvalue support of any edge state, we may assume that $2-2\cos\left(\pi r/n\right)$ is a non-zero eigenvalue that is not in the eigenvalue support of $e_k-e_{k+1}$ for some integer $1\leq r\leq n-1$.  Let $v_r$ denote the eigenvector of $L(P_n)$ such that \[v_rv_r^T=E_r.\] 

Assume that $2-2\cos\left(\pi r/n\right)$ is not in the eigenvalue support of $(k,k+1)$, which means that \[
E_r\left( e_k -e_{k+1}\right)=0.\] Then we know that \[
v_r e_k^T =v_r e_{k+1}^T.\] By Lemma~\ref{pathevt}, we must have that
\[\cos((2k-1)\frac{r\pi}{2n})=\cos\left((2k+1)\frac{r\pi}{2n}\right).\] Using the trigonometric identity
\[
\cos(x)-\cos(y)=-2\sin\left(\frac{x+y}{2}\right)\sin\left(\frac{x-y}{2}\right),\] we know that $r$ must satisfy
\[\cos((2k-1)\frac{r\pi}{2n})-\cos\left((2k+1)\frac{r\pi}{2n}\right)=
-2\sin\left(4k\frac{r\pi}{4n}\right)\sin\left(\frac{r\pi}{2n}\right)=0.\] Thus, we know that either $
kr/n$ or $r/2n$ is an integer. But $1\leq r\leq n-1$ and so \[\frac{kr}{n}=z\] for some positive integer $z$. 

Since $1\leq r\leq n-1$, we know that $z$ must satisfy that \[
1\leq \frac{n}{k}z\leq n-1.\] The number of values of $z$ satisfying the inequality above is the number of non-zero eigenvalues not in the eigenvalue support of $e_k-e_{k+1}$. By Lemma~\ref{pathsymmetric eigen_supp}, we only need to consider the cases when $k=1,2,\cdots,\lfloor\frac{n}{2}\rfloor.$ Since the number of valid $z$ increases as the value of $k$ increases and when $k=\lfloor \frac{n}{2}\rfloor$, the values that $z$ can take is at most \[ 
\lfloor \frac{n}{2}\rfloor-1 \leq \frac{n}{2}-1.\] As stated before, zero is never in the eigenvalue support of a pair state and so, we can conclude that there are at most $n/2$ eigenvalues that are not in the eigenvalue support of an edge state in $P_n$. Therefore, the size of the eigenvalue support of an edge state is at least ${n}/{2}$.\qed

Since we only consider the case when perfect pair state transfer between pair states that at least one of them is an edge state, we conclude the following theorem.

\begin{theorem}
\label{pathpst}
A path graph on $n$ vertices has perfect pair state transfer if and only if  $n=3,4$.
\end{theorem}
\proof By Lemma~\ref{pathevt}, we know that $P_n$ has Laplacian eigenvalue \[0\leq 2-2\cos\frac{\pi r}{n}\leq 4\] for $r = 0,1,\cdots,n-1$. By Corollary~\ref{pedeigendiff}, we know that if an edge state of $P_n$ is periodic, then its eigenvalue support has size at most four. Lemma~\ref{path eigen_supp size} tells us that the eigenvalue support of an edge state of $P_n$ is at least $
{n}/{2}$. Thus, we know that for $n\geq 9$, there is no periodic edge states in $P_n$, which implies that there is no perfect pair state transfer in $P_n$ when $n\geq 9$. Thus, we only need to consider the cases when $n=3,4,5,6,7,8$.

Using Theorem~\ref{pedev} we find that when $n=5,7,8,9$, there is no periodic edge states in $P_n$. Since if an edge state has perfect state transfer, then it must be periodic, which tells us that when $n=5,7,8,9$, there is no perfect pair state transfer in $P_n$.

By computing \[
\abs{\frac{1}{2}(e_a - e_b)^TU(t)(e_c - e_d)} ^2\] for all different vertex-pairs $(a,b),(c,d)$ in $P_3,P_4$ and $P_6$, we find that there is perfect state transfer in $P_3$ and $P_4$. Therefore, there is perfect state transfer in $P_n$ if and only if $n=3,4$.\qed

When $n=3$, there is perfect state transfer between its edges in $P_3$ at time $\pi/2$. When $n=4$, perfect state transfer occurs between two edges on its ends in $P_4$ at time $\sqrt{2}\pi/2$.
\subsection{Comments}
Stevanovi\'c~\cite{Stevanovic2011} and Godsil\cite{Godsil2012} prove that
$P_n$ admits perfect vertex state transfer relative to adjacency matrices if and only if $n=2$ or $3$. Perfect vertex state transfer in $P_2$ happens between its two vertices at time $\pi/2$ and perfect vertex state transfer in $P_3$ happens between its end-vertices at time $\sqrt{2}\pi/2$. 

We proved that $P_n$ admits perfect pair state transfer only when $n=3$ or $4$ and 
\begin{enumerate}[label=(\roman*)]
\item there is perfect state transfer between its edges in $P_3$ at time $\pi/2$,
\item when $n=4$, perfect state transfer occurs between two edges on its ends in $P_4$ at time $\sqrt{2}\pi/2$.
\end{enumerate}
Later in Section~\ref{unsigned}, we will prove an analogous result for quantum walks relative to the unsigned Laplacians in paths with initial states of the form $e_a+e_b$. That is, $P_3,P_4$ are the only paths where perfect state transfer relative to the unsigned Laplacians occurs and it occurs between the end-edges of $P_3,P_4$ at time $\pi/2$, $\sqrt{2}\pi/2$ respectively .

Notice also that $P_2, P_3$ are the line graphs of $P_3,P_4$ respectively. In $P_3$ and its line graph $P_2$, perfect state transfer always occurs at the same time $\pi/2$ between the same pair of edges and their corresponding pair of vertices in the line graph. This happens regardless of our choice of Hamiltonian or form of the initial state. We can make the same observations about perfect state transfer in $P_4$ ant its line graph $P_3$.

We know that $C_4$ is the only cycle that admits perfect state transfer relative to adjacency matrices, Laplacians. We will show in next section that $C_4$ is the only cycle that admits perfect state transfer relative to the unsigned Laplacians, where the initial state is in the plus state form. No matter our choice of Hamiltonians and form of initial state, perfect state transfer in $C_4$ happens at the same time $\pi/2$ between pairs of opposite edges or vertices.

Let $G$ be a regular graph with valency $k$, then the Laplacian matrix of $G$ is \[
L=kI-A.\] Then the transition matrix for pair state transfer is 
\[U(t) = \exp\left( it (kI-A)\right) = e^{itk} e^{-itA}.\] A similar argument works for the unsigned Laplacians. Thus we can conclude that the continuous quantum walks generated by the adjacency matrices, the Laplacians and the unsigned Laplacians are equivalent up to a phase factor. 

Although the transition matrices with respect to the adjacency matrices, the Laplacians and the unsigned Laplacians of cycles are equivalent, it is still surprising that different forms of initial states actually do not affect perfect state transfer in $C_4$. 
Also, notice that $C_4$ is the line graph of itself and there is a correspondence between pairs of PST-edges and pairs of PST-vertices.

It may seem that there is a correspondence between perfect edge state transfer in a graph and perfect vertex state transfer in its line graph. However, that is not true for most graphs. So far, paths and cycles are the only examples we have found where the correspondence can be observed.

\section{Unsigned Laplacian}
\label{unsigned}
Let $G$ be a graph. The unsigned Laplacian of $G$ is matrix $L_+(G)$ such that \[
L_+(G)=\Delta(G) +A(G).\] When we use $L_+(G)$ as Hamiltonian in a quantum walk,  the pair of vertices $(a,b)$ of $G$ is associated with the state \[
e_a+e_b,\] which we call ``plus state".

Every time we refer to plus states, we use the unsigned Laplacian of a graph as Hamiltonian unless stated explicitly otherwise. We define analogously that there is perfect plus state transfer between $e_a+e_b$ and $e_c+e_d$ if and only if \[
U(t)(e_a+e_b)=\exp\left(itL_+\right)(e_a+e_b)=\gamma (e_c+e_d),\] for some complex constant $\gamma$ with norm $1$. Also, a plus state $e_a+e_b$ is periodic if and only if it has perfect plus state transfer to itself at some time $t$.

Since the main case of interest in this paper is the case when the Laplacian of a graph is used as Hamiltonian, it is natural to question if there will be perfect state transfer between a pair state and a plus state when we use the Laplacian as Hamiltonian. The answer is no.
\begin{theorem} Let $G$ be a graph with $a,b,c,d\in V(G)$.
There is no perfect state transfer between a state of the form $e_a+e_b$ and a state of the form $e_c-e_d$ in $G$ when the Laplacian of $G$ is used as Hamiltonian of the quantum walk.
\end{theorem}
\proof
We know that $0$ will always be a eigenvalue of the Laplacian of $G$ with the all-ones vector being its eigenvector. Thus, we know $0$ will never be in the eigenvalue support of $e_a-e_b$ while $0$ is always in the eigenvalue support of $e_c+e_d$. It follows that $e_a-e_b$ and $e_c+e_d$ do not have the same eigenvalue support, which implies that they are not strongly cospectral. By Theorem~\ref{pst-strongly cospectral}, we can conclude that there is no perfect state transfer between a state of the form $(e_a-e_b)$ and a state of the form $(e_c+e_d)$ using Laplacian as Hamiltonian.\qed

\begin{table}[H]
\begin{center}
\scalebox{0.85}{
\begin{tabular}{|c|c|c|c|c|c|}
\hline
$G_n$& total  & Lap. PED&Lap. PST
& Unsigned PED  & Unsigned PST \\
\hline
$G_5$ & $21$ & $18\text{ }(85.7\%)$&$6\text{ }(28.6\%)$ & $4\text{ }(19.0\%)$&$0\text{ }(0\%)$\\
\hline
$G_6$ & $112$ & $86\text{ }(76.8\%)$&$25\text{ }(22.3\%)$ & $21\text{ }(18.8\%)$&$4\text{ }(3.6\%)$\\
\hline
$G_7$ & $853$ & $513\text{ }(60.1\%)$&$94\text{ }(11.0\%)$ & $23\text{ }(2.7\%)$&$2\text{ }(0.2\%)$\\
\hline
$G_8$ & $11117$ & $5164\text{ }(46.5\%)$   & $673\text{ }(6.0\%)$ & $55\text{ }(0.5\%)$ &  $14\text{ }(0.1\%)$\\
\hline
\end{tabular}}
\caption{the Number of Graphs with PST and Periodic States}
\label{sign and unsigned pst ped}
\end{center}
\end{table}

Despite the huge gap between the number of PST pairs in terms of pair state transfer and plus state transfer, when the underlying graph is a bipartite graph, perfect state transfer in terms of pair states and plus states are equivalent.
\begin{lemma}
\label{bipartiteLL+}
Let $G$ be a bipartite graph with two parts $B_1,B_2$ and $A,\Delta$ denote the adjacency matrix and the degree matrix of $G$ respectively.  Let $D$ be block matrix such that \[D=\begin{pmatrix}
-I & \mathbf{0}\\
\mathbf{0} & I
\end{pmatrix}\] indexed by the vertices of $B_1,B_2$ in the order . Then we have \[
D(\Delta-A) D=\Delta+A.\qed\]
\end{lemma}

\begin{theorem}
\label{bipartite pst}
Let $G$ be a bipartite graph with parts $B_1, B_2$ and vertices $a,c\in B_1$ and $b,d\in B_2$. There is perfect pair state transfer between $(e_a-e_b)$ and $(e_c-e_d)$ if and only if there is perfect plus state transfer between $(e_a+e_b)$ and $(e_c+e_d)$.
\end{theorem}
\proof
Let $\Delta$ denote the degree matrix of $G$ and $A$ denote the adjacency matrix of $G$. From the Lemma~\ref{bipartiteLL+}, we know that\[
D(\Delta-A)D=\Delta+A\] and inserting $DD=I$ between $m$ copies of $\Delta-A$, we have
 \begin{align*}
D(\Delta-A)^m D &= D(\Delta-A)DD(\Delta-A)DD\cdots (\Delta-A)D\\[1em]
& = (\Delta+A)^m
\end{align*} for any non-negative integer $m$. Then we see that \begin{align*}
DU(t)D =D\exp(itL)D&=D \sum_{m=0}^\infty \left(\frac{(it)^m}{m!}\left(\Delta-A\right)^m\right) D\\[1em]
&=\sum_{m=0}^\infty \frac{(it)^m}{m!}D\left(\Delta-A\right)^m D\\[1em]
&=\sum_{m=0}^\infty \frac{(it)^m}{m!}\left(\Delta+A\right)^m\\[1em]
&=\exp(itL_+)
\end{align*}
Note that since $a,c\in B_1$ and $b,d\in B_2$, we have that \[
D(e_a-e_b)=-( e_a+e_b),\quad D(e_c-e_d)=-( e_c+e_d).\] There is perfect state transfer between $(e_a-e_b)$ and $(e_c-e_d)$ using Laplacian if and only if there exist $\tau$ such that 
\[U(\tau)(e_a-e_b) = \gamma (e_c-e_d)\] for some $\abs{\gamma}=1.$ Applying $D$ on both sides of the equation above, we have that 
\[
DU(\tau)(e_a-e_b) =D\left( \gamma (e_c-e_d)\right).\]
Again using $DD=I$, we can rewrite the equation above as \[
DU(\tau)DD(e_a-e_b) = \gamma D (e_c-e_d).\] This gives us that \begin{gather*}
-\exp(i\tau L_+)(e_a+e_b)=-\gamma(e_c+e_d),\\[1em]
\exp(i\tau L_+)(e_a+e_b)=\gamma(e_c+e_d),
\end{gather*}which is equivalent to perfect plus state transfer between $(e_a+e_b)$ and $(e_c+e_d)$ using unsigned Laplacian. This completes our proof.\qed

Next we discuss perfect plus state transfer in cycles and paths. Like the case in pair state transfer, we exclude the perfect state transfer between both non-edge state.

Since we proved that $C_4$ is the only cycle where perfect pair state transfer occurs, one immediate result from the theorem above is that there is no perfect plus state transfer in $C_n$ when $n$ is even and the only exception is when $n=4$. In $C_4$, there is perfect plus state transfer between opposite edges at time $\frac{\pi}{2}$.

A plus-state analogue of Corollary~\ref{noPST in odd edge-transitive} implies that there is no perfect plus state transfer in $C_n$ when $n$ is odd. Therefore, we get the following theorem.

\begin{theorem}
There is perfect plus state transfer in $C_n$ if and only if $n=4$. \qed
\end{theorem}

Since $P_3,P_4$ are the only paths where perfect pair state transfer occurs, Theorem~\ref{bipartite pst} gives the result below.
\begin{theorem}
\label{unsigned path pst}
A path graph on $n$ vertices has perfect plus state transfer if and only if  $n=3,4$.\qed
\end{theorem}

Perfect pair  state transfer and perfect plus state transfer happen between the same pairs of edges at the same time in $P_3$ and $P_4$.  When $n=3$, there is perfect plus state transfer between its edges in $P_3$ at time $\pi/2$. When $n=4$, perfect plus state transfer occurs between two edges on its ends in $P_4$ at time $\sqrt{2}\pi/2$.

\section{Open Questions}
Checking all the trees on up to $16$ vertices, there is no perfect pair state transfer and there are only four types of graphs that contain periodic pair states:\begin{enumerate}
\item Star graphs $K_{1,n}$;
\item Double stars;
\item Paths;
\item The figure below.
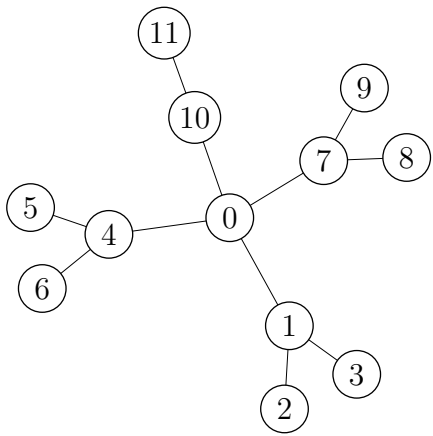
\begin{figure}[H]
\begin{center}
\begin{tikzpicture}
\definecolor{cv0}{rgb}{0.0,0.0,0.0}
\definecolor{cfv0}{rgb}{1.0,1.0,1.0}
\definecolor{clv0}{rgb}{0.0,0.0,0.0}
\definecolor{cv1}{rgb}{0.0,0.0,0.0}
\definecolor{cfv1}{rgb}{1.0,1.0,1.0}
\definecolor{clv1}{rgb}{0.0,0.0,0.0}
\definecolor{cv2}{rgb}{0.0,0.0,0.0}
\definecolor{cfv2}{rgb}{1.0,1.0,1.0}
\definecolor{clv2}{rgb}{0.0,0.0,0.0}
\definecolor{cv3}{rgb}{0.0,0.0,0.0}
\definecolor{cfv3}{rgb}{1.0,1.0,1.0}
\definecolor{clv3}{rgb}{0.0,0.0,0.0}
\definecolor{cv4}{rgb}{0.0,0.0,0.0}
\definecolor{cfv4}{rgb}{1.0,1.0,1.0}
\definecolor{clv4}{rgb}{0.0,0.0,0.0}
\definecolor{cv5}{rgb}{0.0,0.0,0.0}
\definecolor{cfv5}{rgb}{1.0,1.0,1.0}
\definecolor{clv5}{rgb}{0.0,0.0,0.0}
\definecolor{cv6}{rgb}{0.0,0.0,0.0}
\definecolor{cfv6}{rgb}{1.0,1.0,1.0}
\definecolor{clv6}{rgb}{0.0,0.0,0.0}
\definecolor{cv7}{rgb}{0.0,0.0,0.0}
\definecolor{cfv7}{rgb}{1.0,1.0,1.0}
\definecolor{clv7}{rgb}{0.0,0.0,0.0}
\definecolor{cv8}{rgb}{0.0,0.0,0.0}
\definecolor{cfv8}{rgb}{1.0,1.0,1.0}
\definecolor{clv8}{rgb}{0.0,0.0,0.0}
\definecolor{cv9}{rgb}{0.0,0.0,0.0}
\definecolor{cfv9}{rgb}{1.0,1.0,1.0}
\definecolor{clv9}{rgb}{0.0,0.0,0.0}
\definecolor{cv10}{rgb}{0.0,0.0,0.0}
\definecolor{cfv10}{rgb}{1.0,1.0,1.0}
\definecolor{clv10}{rgb}{0.0,0.0,0.0}
\definecolor{cv11}{rgb}{0.0,0.0,0.0}
\definecolor{cfv11}{rgb}{1.0,1.0,1.0}
\definecolor{clv11}{rgb}{0.0,0.0,0.0}
\definecolor{cv0v1}{rgb}{0.0,0.0,0.0}
\definecolor{cv0v4}{rgb}{0.0,0.0,0.0}
\definecolor{cv0v7}{rgb}{0.0,0.0,0.0}
\definecolor{cv0v10}{rgb}{0.0,0.0,0.0}
\definecolor{cv1v2}{rgb}{0.0,0.0,0.0}
\definecolor{cv1v3}{rgb}{0.0,0.0,0.0}
\definecolor{cv4v5}{rgb}{0.0,0.0,0.0}
\definecolor{cv4v6}{rgb}{0.0,0.0,0.0}
\definecolor{cv7v8}{rgb}{0.0,0.0,0.0}
\definecolor{cv7v9}{rgb}{0.0,0.0,0.0}
\definecolor{cv10v11}{rgb}{0.0,0.0,0.0}
\Vertex[style={minimum size=1.0cm,draw=cv0,fill=cfv0,text=clv0,shape=circle},LabelOut=false,L=\hbox{$0$},x=2.6482cm,y=2.5439cm]{v0}
\Vertex[style={minimum size=1.0cm,draw=cv1,fill=cfv1,text=clv1,shape=circle},LabelOut=false,L=\hbox{$1$},x=3.4404cm,y=1.1062cm]{v1}
\Vertex[style={minimum size=1.0cm,draw=cv2,fill=cfv2,text=clv2,shape=circle},LabelOut=false,L=\hbox{$2$},x=3.3757cm,y=0.0cm]{v2}
\Vertex[style={minimum size=1.0cm,draw=cv3,fill=cfv3,text=clv3,shape=circle},LabelOut=false,L=\hbox{$3$},x=4.3361cm,y=0.4609cm]{v3}
\Vertex[style={minimum size=1.0cm,draw=cv4,fill=cfv4,text=clv4,shape=circle},LabelOut=false,L=\hbox{$4$},x=1.0414cm,y=2.3227cm]{v4}
\Vertex[style={minimum size=1.0cm,draw=cv5,fill=cfv5,text=clv5,shape=circle},LabelOut=false,L=\hbox{$5$},x=0.0cm,y=2.6794cm]{v5}
\Vertex[style={minimum size=1.0cm,draw=cv6,fill=cfv6,text=clv6,shape=circle},LabelOut=false,L=\hbox{$6$},x=0.1552cm,y=1.6026cm]{v6}
\Vertex[style={minimum size=1.0cm,draw=cv7,fill=cfv7,text=clv7,shape=circle},LabelOut=false,L=\hbox{$7$},x=3.897cm,y=3.303cm]{v7}
\Vertex[style={minimum size=1.0cm,draw=cv8,fill=cfv8,text=clv8,shape=circle},LabelOut=false,L=\hbox{$8$},x=5.0cm,y=3.3455cm]{v8}
\Vertex[style={minimum size=1.0cm,draw=cv9,fill=cfv9,text=clv9,shape=circle},LabelOut=false,L=\hbox{$9$},x=4.4382cm,y=4.2693cm]{v9}
\Vertex[style={minimum size=1.0cm,draw=cv10,fill=cfv10,text=clv10,shape=circle},LabelOut=false,L=\hbox{$10$},x=2.182cm,y=3.8783cm]{v10}
\Vertex[style={minimum size=1.0cm,draw=cv11,fill=cfv11,text=clv11,shape=circle},LabelOut=false,L=\hbox{$11$},x=1.7769cm,y=5.0cm]{v11}
\Edge[lw=0.01cm,style={color=cv0v1,},](v0)(v1)
\Edge[lw=0.01cm,style={color=cv0v4,},](v0)(v4)
\Edge[lw=0.01cm,style={color=cv0v7,},](v0)(v7)
\Edge[lw=0.01cm,style={color=cv0v10,},](v0)(v10)
\Edge[lw=0.01cm,style={color=cv1v2,},](v1)(v2)
\Edge[lw=0.01cm,style={color=cv1v3,},](v1)(v3)
\Edge[lw=0.01cm,style={color=cv4v5,},](v4)(v5)
\Edge[lw=0.01cm,style={color=cv4v6,},](v4)(v6)
\Edge[lw=0.01cm,style={color=cv7v8,},](v7)(v8)
\Edge[lw=0.01cm,style={color=cv7v9,},](v7)(v9)
\Edge[lw=0.01cm,style={color=cv10v11,},](v10)(v11)
\end{tikzpicture}
\caption{$(0,10)$ is the only periodic pair in this graph}
\end{center}
\end{figure}
\end{enumerate}

Coutinho and Liu in \cite{Coutinho2014} proved that there is no perfect vertex state transfer in trees with more than two vertices when the Laplacian is the Hamiltonian. We suspect a similar result holds 
for pair state transfer.
(Our personal feeling is that there is no perfect pair state transfer on trees, but we have not found a way to prove this.)

Another question we would like to answer is that how different Hamiltonians and different initial states affect state transfer.
\begin{table}[H]
\begin{center}
\scalebox{0.8}{
\begin{tabular}{|c|c|c|c|c|c|c|c|}
\hline
$G_n$ & Total  & \vtop{\hbox{\strut$A(G)$ with}\hbox{\strut vertex states}}& Prop. & \vtop{\hbox{\strut$L(G)$ with}\hbox{\strut edge states}} & Prop.& \vtop{\hbox{\strut$L_+(G)$ with}\hbox{\strut plus states}}& Prop.\\
\hline
$G_5$ & $21$ &$1$& $4.8\%$ &$6$ & $28.6\%$ & $0$ & $0\%$ \\
\hline
$G_6$ & $112$ & $1$& $0.9\%$ &$27$& $24.1\%$ & $4$ & $3.6\%$\\
\hline
$G_7$ & $853$ & $1$ & $0.1\%$ & $104$  &  $12.2\%$ & $2$ & $0.2\%$ \\
\hline
$G_8$ & $11117$ & $5$ & $0.004\%$ &$779$  & $7.0\%$& $14$ & $0.1\%$\\
\hline
\end{tabular}}
\caption{the Number of graphs with PST in Different Settings}
\label{diff ped state}
\end{center}
\end{table}
Table~\ref{diff ped state} shows that the number of graphs with adjacency vertex-state PST, Laplacian pair-state PST and unsigned Laplacian plus-state PST followed with the corresponding proportions from left to right in order.

We can see that on a set of graphs, different choices of Hamiltonian and different forms of initial state strongly affect the number of graphs that have perfect state transfer. However, as shown in Section~\ref{specialcases} and Section~\ref{unsigned}, when the underlying graphs are bipartite graphs and odd cycles, perfect edge state transfer and perfect plus state transfer are equivalent. We can see that perfect state transfer in certain classes of graphs is invariant under different Hamiltonians and initial states.

We want to answer the question that given a specific Hamiltonian for a graph $G$, which form of the initial states (i.e.,  vertex states, pair state, plus states) gives us the most perfect state transfer pairs in $G$. On the other hand, given a specific initial state of a graph $G$, we would like to know which Hamiltonian, (i.e. adjacency matrix of $G$, Laplacian of $G$, unsigned Laplacian of $G$) has the advantage of producing the most perfect state transfer pairs in $G$.  Also, we would like to know that  besides bipartite graphs and odd cycles, if there is any other classes of graphs such that different choices of Hamiltonians and initial states do not affect on PST pairs on graphs.

\bibliographystyle{plain}

\end{document}